\documentclass[twoside,a4paper,titlepage]{article}

\usepackage[noend]{algpseudocode}
\usepackage{enumitem}
\usepackage{graphicx}
\usepackage{textcomp}
\usepackage{makeidx}
\usepackage{wrapfig}
\usepackage{amsmath}
\usepackage{amssymb}
\usepackage{amsthm}
\usepackage{bbm}

\renewcommand{\mathbb}[1]{\mathbbm{#1}}

\algrenewcommand\alglinenumber[1]{
        {\sf\footnotesize#1}}
\algblockdefx[When]{When}{EndWhen}[1]{\textbf{when} #1}{}
\algblockdefx[Initially]{Initially}{EndInitially}[0]{\textbf{initially}}{}

\newcommand{\N}{\mathbb{N}}
\newcommand{\W}{\mathcal{W}}
\newcommand{\K}{\mathcal{K}}
\newcommand{\keyword}[1]{\textbf{#1}}

\newcommand{\halts}{\mathclose\downarrow}
\newcommand{\divs}{\mathclose\uparrow}
\newcommand{\eqT}{\equiv_{\mathrm{T}}}
\newcommand{\leqT}{\leq_{\mathrm{T}}}
\newcommand{\leqm}{\leq_{\mathrm{m}}}
\newcommand{\leqwtt}{\leq_{\mathrm{wtt}}}
\newcommand{\leT}{<_{\mathrm{T}}}

\newcommand{\ind}{\mathbb{1}}
\newcommand{\mathsc}[1]{\ensuremath{\text{\textsc{#1}}}}
\newcommand{\Pfin}{\mathfrak{P}_{\text{fin}}}
\newcommand{\PP}{\mathfrak{P}}
\newcommand{\dom}{\operatorname{dom}}
\newcommand{\im}{\operatorname{Im}}

\newcounter{main}
\theoremstyle{plain}

\newtheorem{prop}[main]{Proposition}
\newtheorem{thm}[main]{Theorem}
\newtheorem{cor}[main]{Corollary}
\theoremstyle{definition}
\newtheorem{dfn}[main]{Definition}
\newtheorem{exa}[main]{Example}
\newtheorem{rem}[main]{Remark}
\newtheorem{prob}[main]{Problem}
\newtheorem{oprob}[main]{Open Problem}

\author{Bas Westerbaan \\ \texttt{bas@westerbaan.name}}
\title{On effective undecidability and Post's problem}

\makeindex

\begin{document}

\maketitle

\section{Introduction}
A subset~$A \subseteq \N$ of the natural numbers is called
\keyword{decidable}\index{decidable},
if its characteristic function~$\ind_A$, defined by
\[ \ind_A (n) = \begin{cases}
              1 & n \in A \\
              0 & n \notin A,
             \end{cases} \]
is computable. 
\keyword{Computable}\index{computable|see{decidable}},
\keyword{recursive}\index{recursive|see{decidable}}
and \keyword{effective}\index{recursive|see{decidable}} are
synonyms.

A subset~$A \subseteq \N$ is
called~\keyword{recursively enumerable}\index{enumerable!recursively}
\index{r.e.|see{recursively enumerable}}
if it is finite
or there exists a computable~$a\colon \N \to \N$ that enumerates it. That is:
\[ A = \{ a(0), \, a(1),\, \ldots \} .\]
Similarly
\keyword{computably enumerable}
\index{enumerable!computably|see{recursively}} and
\keyword{effectively enumerable}
\index{enumerable!effectively|see{recursively}}
are synonyms.
Every decidable set is recursively enumerable, but there are a lot of
recursively enumerable sets that are not decidable.

An important example of such a set is~$\K$.%
\footnote{$\K$ appears (implicitely) in
the work of G\"odel \cite{god31}, Turing \cite{tur36}
and Kleene \cite{klee36}}.
Let~$\varphi^1_0,\, \varphi^1_1,\, \ldots$ be a standard enumeration of all
computable partial functions from~$\N$ to~$\N$ and
let~$\varphi^1_{37}(n)\halts$
denote that~$\varphi^1_{37}$ is defined on~$n$.
Then we define\index{K@$\K$}
\[ \K = \{ e; \, e \in \N; \,\varphi^1_e(e) \halts \}. \]

We can find a recursive~$r\colon \N \to \N$ such
that
\begin{equation*}
\varphi^1_{r(e)}(x) \halts \iff \varphi^1_e(x)=0.
\end{equation*}
Then
\begin{equation*}
r(e) \in \K \iff \varphi^1_{r(e)}(r(e)) \halts
        \iff \varphi^1_e(r(e)) = 0
\end{equation*}
Hence~$\varphi^1_e(r(e)) \neq \ind_\K(r(e))$, for any~$e$. Thus
no~$\varphi^1_e$ can compute~$\K$.

Two famously undecidable sets are
\begin{align*}
\mathsc{pa} & = \{ \# \psi; \,\psi \in L_{\mathsc{N}}; \,
                        \text{PA} \vdash \psi \} \\
\mathsc{n} & = \{ \# \psi; \,\psi \in L_{\mathsc{N}}; \,
                        \left< \N, +, \cdot, 0, 1, \leq \right> \vDash \psi \},
\end{align*}
where~$L_{\mathsc{N}}$ is the set of first-order formulas about the natural
numbers and~$\# \psi$ is the G\"odelcode of~$\psi$.
Thus~$\mathsc{pa}$ is the set of natural numbers coding a sentence that is 
provable from the axioms of Peano and~$\mathsc{n}$ is the (bigger) set of
natural numbers coding a sentence that is  true about the natural numbers.

The undecidability of~$\mathsc{n}$ is a celebrated result by Alan
Turing\cite{tur36}\footnote{It was also discovered independently
by Alonzo Church.\cite{chu36}} and a definitive blow to Gottfried
Leibniz' and David Hilbert's dream that (mathematical) truth might
be established mechanically.

The set~$\K$ is as undecidable as~$\mathsc{pa}$ in the following
sense.  If there were an oracle which would answer every query
of the form ``is $n\in\K$?'', then we can give an algorithm using
this oracle to compute~$\mathsc{pa}$.  And vice versa.

The set~$\mathsc{n}$, however, is ``more undecidable'' than~$\mathsc{pa}$.
That is: even if we allow an algorithm to use an oracle for~$\mathsc{pa}$,
it still cannot compute~$\mathsc{n}$. It is not even recursively
enumerable.

This leads to the study of degrees of undecidability. In 1944 Emil Post
asked the following question:
\begin{prob}[Post \cite{post44}]
Is there a recursively enumerable set that is undecidable but strictly
``less undecidable'' than~$\K$?
\end{prob}
\begin{wrapfigure}{r}{0.3\textwidth}
\begin{center}
\includegraphics[width=0.28\textwidth]{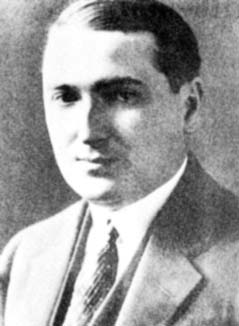}
\end{center}
\caption{Emil Post}
\end{wrapfigure}

This problem was solved 14 years later by Albert Muchnik\cite{much58}
(and independently a year later by Richard Friedberg\cite{frie57})
using a method, now called the \emph{priority method}, which has
become a main tool in recursion theory.

The set they constructed is, in a sense,
artificial.  One would hope for a
set that emerges more naturally, such as~$\mathsc{pa}$.

One uncanny aspect of the set constructed by Friedberg and Muchnik, say~$X$,
is that we only have an indirect method to prove that~$X$ differs
from any given decidable set~$D$.
The method suggests a number~$n_1$ that might be an example of a number
that is in~$D$, but is not in~$X$.  When~$n_1$ nevertheless
becomes a member of~$X$,
the method yields another number~$n_2$, which again might be an example
of a number that is in~$D - X$. When~$n_2$ appears
in~$X$, the method yields a~$n_3$ which might be in~$D-X$. Et cetera.
This, however, does not go on indefinitely.
There is a computable~$f\colon \N \to \N$ such that for the~$e$th
decidable set, the method requires less than~$f(e)$ tries.

In Post's search for a solution and also in later work of Martin
and others it appears that every set that is in some reasonable
sense `effectively undecidable' ends up being as undecidable as~$\K$.

In this thesis, we will look at some known and some previously 
uninvestigated notions of effective undecidability.  We try to discover
how far we can stretch effective undecidability in the hope to get a
more tractable solution to Post's problem, than that of Friedberg
and Muchnik.
\clearpage
\tableofcontents
\clearpage
\section{Basic recursion theory}
In this section, we will briefly review the basics of recursion theory.
The reader familiar with the matter, is suggested to skip to
Subsection~\ref{ss:conventions}.
\subsection{Models of computation and the Church-Turing Thesis}
In the introduction, we informally talked about computable functions.
Can we define these formally?
There have been dozens of wildly varying suggestions for formalization
of the notion of computation.  All suggestions so far are proven
equivalent.\footnote{See Theorem~I7.12 of~\cite{odi1}}
The Church-Turing Thesis asserts that the notion of computation is captured
by any of those.

We will briefly review three formalizations of computation:~$\mu$-recursive,
$\lambda$-definable and~Turing-computable.

\subsubsection{$\mu$-recursive functions}
\index{recursive!$\mu$-}%
Let~$\mathcal{F}_k$ be the set of partial functions from~$\N^k$ to~$\N$
and~$\mathcal{F} = \bigcup_{k\in\N} \mathcal{F}_{k+1}$.
The set of~$\mu$-recursive functions is a subset of~$\mathcal{F}$ defined
inductively as follows:\cite{vel87}
\begin{enumerate}
\item Zero function. $0 \mapsto 0$ is~$\mu$-recursive.
\item Projections. For any~$k$ and~$i \in \{1,\ldots, k\}$
the function~$(n_1, \ldots, n_k) \mapsto n_i$ is~$\mu$-recursive.
\item Successor. The function~$n \mapsto n+1$ is~$\mu$-recursive.
\item Composition. For any~$k$, $l$ and~$\mu$-recursive~$k$-ary~$h$
and~$\mu$-recursive~$l$-ary~$g_1,\ldots,g_k$ the function
\[ (n_1, \ldots, n_l) \mapsto h(g_1(n_1, \ldots, n_l), \ldots,
                g_k(n_1, \ldots, n_l)) \]
is also~$\mu$-recursive.
\item Primitive recursion. 
Given any $\mu$-recursive $(k+2)$-ary function~$g$,
the~$(k+1)$-ary function~$f$ such that
\begin{align*}
        f(0, x_1, \ldots, x_k) & = g(0, 0, x_1, \ldots, x_k) \\
        f(n+1, x_1, \ldots, x_k) & = g(n+1,
                        f(n, x_1, \ldots, x_k), x_1, \ldots, x_k)
\end{align*}
is also~$\mu$-recursive.
\end{enumerate}

There is one rule left: minimization. The partial functions
defined by only the first 5 rules are in fact total and 
called the primitive recursive functions.\index{recursive!primitive}
In the early 1900s it was believed that all total computable functions
can be defined by primitive recursion.\cite{dotz91}
In 1923, Wilhelm Ackermann published a counterexample\footnote{
This is Ackermann's function. It is total and computable, but
not primitive recursive.
\[ A(m,n) = \begin{cases}
            n + 1 & \text{if $m=0$} \\
            A(m-1,1) & \text{if $m>0$ and $n=0$} \\
            A(m-1,A(m,n-1)) & \text{if $m>0$ and $n>0$}
            \end{cases} \]}.
To capture all computable functions, Kleene suggested the final rule
for~$\mu$-recursive functions:
\begin{enumerate}[resume]
\item Minimisation. Given any~$\mu$-recursive $(k+1)$-ary function~$g$,
the~$k$-ary function~$f$ such that
\begin{align*}
 f(x_1, \ldots, x_k) = z \iff &
        f(z, x_1, \ldots, x_k) = 0 \text { and }\\
        &\forall n < z \exists y \neq 0 [ f(n, x_1, \ldots, x_k) = y ]
\end{align*}
is also~$\mu$-recursive.  In general~$f$ is partial.
\end{enumerate}

Although the definition of~$\mu$-recursive functions is elegant,
it is not intuitively clear that they are exactly the computable functions:
we are a priori not sure that there is no counterexample like that
Ackermann for the primitive recursive functions.
However, rest assured: the other two equivalent formalizations that follow
are rather more convincing in capturing the notion of computability.

\subsubsection{Lambda Calculus}
\index{lambda-definable@$\lambda$-definable}%
Another formalization is the Lambda Calculus.\cite{bar00} Alonzo Church proved
in 1936 that there is no~$\lambda$-definable function that decides
whether two given~$\lambda$-terms are $\beta$-equal.\cite{chu36}
These two notions will be defined later on.  If~$\lambda$-definable
is the same as computable, this proves that there are classes of mathematical
statements that cannot be solved algorithmically.

Given a function~$f$ and an argument~$x$, creating~$f(x)$ is called
\emph{(function) application}. Given some expression with a free
variable (for instance~$y^2 - \int_0^1 x^2 \,\textrm{d}x$) creating a
function in that free variable (in our example: $f(y)=y^2 - \int_0^1
x^2 \,\textrm{d}x$) is called \emph{(function) abstraction}.

Usually, when one considers some kind of functions, their arguments
are of a different type than the functions themselves. Linear functions
act on vectors, functors act on structures, et cetera.

The~$\lambda$-calculus is a language to describe functions
that act on themselves and (because that on itself would
not be very interesting) allow (function) abstractions.

The \emph{words} in the calculus are the \emph{lambda terms}.
\begin{itemize}
\item $x$ is a variable symbol.
\item If~$v$ is a variable symbol, then~$v'$ is a variable symbol.
\item If~$v$ is a variable symbol, then~$v$ is a lambda term.
\item If~$M$ and~$N$ are lambda terms, then the word~$(MN)$ is a lambda term.
        This corresponds to function application:~$M$ applied to~$N$.
\item If~$v$ is a variable symbol and~$M$ is a lambda term, then
                the word $(\lambda v.\, M)$ is a lambda term.
        This corresponds to function abstraction:~$v \mapsto M$.
\end{itemize}

Thus~$(\lambda x.\, (xx))$ corresponds to the function that applies its
argument to itself. By convention, we will write~$x,y,z,\ldots$
for variables; leave out redundant parenthesis; let application
bind stronger than abstraction; associate application to the left
and abstraction to the right \emph{and} write~$\lambda a_1 \cdots a_n.\, M$
as shorthand for~$\lambda a_1.\, \cdots \lambda\, a_n. M$.
If two terms are equal up to renaming of variables that respects
binding of~$\lambda$, we consider them equal.  Example:
\begin{align*}
 \lambda xy.\, y(xx) & = (\lambda x. \, (\lambda x'.\, (x' (x x)))) \\
                     & = \lambda yx.\, x(yy) \\
                     & \neq \lambda yx.\, y(xx).
\end{align*}
To capture the intention to model function, we study the following reduction:
\[ (\lambda v.\, M )N \rightarrow_\beta M[v := N]. \]
Here~$M[v:=N]$ is understood to mean ``$M$ with all unbound occurances of~$v$
are replaced by~$N$'' if~$M$ and~$N$ do not share variables.  If~$M$
and~$N$ do share variables, we can always find equivalent terms
by renaming variables such that this is not the case.
Also if for some lambda term~$M$ there is a subterm~$N$
such that~$N \rightarrow_\beta N'$ then~$M \rightarrow_\beta M'$,
where~$M'$ is~$M$ with that subterm~$N$ replaced by~$N'$.
We write~$M \rightarrow_\beta^* N$ for the
reflexive transitive closure of~$\rightarrow_\beta$ and~$=_\beta$ for the
symmetric transitive closure.
Example:
\begin{align*}
& (\lambda x.\,xx) ((\lambda xyz.\, x(yz)) x (\lambda x. \, xx)) \\
& \rightarrow_\beta  (\lambda xyz.\, x(yz))x (\lambda x.\,xx) ((\lambda xyz.\,
                x(yz)) x (\lambda x. \, xx)) \\
& \rightarrow_\beta  x ((\lambda x.\,xx) ((\lambda xyz.\,
                x(yz)) x (\lambda x. \, xx))) \\
& \rightarrow_\beta  x((\lambda xyz.\, x(yz))x (\lambda x.\,xx) ((\lambda xyz.\,
                x(yz)) x (\lambda x. \, xx))) \\
& \leftarrow_\beta x ((\lambda x.\,xx) ((\lambda xyz.\, x(yz)) x
        (\lambda x. \, xx))).
\end{align*}
We can elegantly represent the natural numbers and their operation in
the lambda calculus using so-called Church-numerals:
$\underline{0} = \lambda so.\, o$ and  $\underline{n+1}  = 
                \lambda so.\,\underline{n} s(so)$.
Thus~$\underline 3 =_\beta \lambda so.\, s(s(so))$.
Let~$\oplus = \lambda nmso. \, n s (m s o)$
and~$\otimes = \lambda nmso. \, n (\lambda c.\, msc) o$.
It is not hard to verify that~$\underline {n + m} 
=_\beta \oplus \underline{n} \underline{m}$
and~$\underline {n \cdot m} 
=_\beta \otimes \underline{n} \underline{m}$ for all~$n,m\in\N$.

A partial function~$f\colon \N \to \N$ is
called~$\lambda$-definable\cite{chu33}\cite{klee35} if there
is a lambda term~$M_f$ such that for all~$n$ and~$m$:
\[ f(n) = m \iff M_f \underline {n} =_\beta \underline{m} .\]

\subsubsection{Turing Machines}
\index{computable!Turing}%
The~$\lambda$-calculus has found many applications especially in theoretical
computerscience. Also there are many programming languages modelled on
the calculus. These are called \emph{functional languages}.
However, at the time (and still on the present day), many
consider~$\lambda$-calculus to be a weird exercise.

A year after Chuch, Alan Turing independently proved that
mathematical truth cannot be established algorithmically.
To do this, he introduced Turing Machines. 
When the famous logician Kurt G\"odel heard of Turing's work, he
promptly stated that Turing finally convincingly defined computability
with his machines.
Computer processors and
imperative programming languages are based on the Turing Machine and
are vastly more popular than their functional counterparts.

A Turing Machine is an idealized computer.
It consists of a tape and a read/write-head on that tape.
The tape extends infinitely in both directions and consists out of blocks.
Each block can hold the value~$0$ or~$1$.  Thus, the state of the tape
can be modeled by a map~$\mathbb{Z} \to \{0,1\}$.
The read/write-head holds a program.  A program consists of a finite
number of states. For each state~$s$ it is specified which action to take
when the head is on a block with a~$0$ and when the head is on a block
with a~$1$. An action consists of either:
\begin{enumerate}
\item        
writing a~$0$ or a~$1$; then moving the head one block to the left
or to the right \emph{and} then changing the state of the machine \emph{or}
\item
stopping the machine.
\end{enumerate}
To \emph{run the machine on input~$n$} for~$n\in\N$
we start with a tape consisting of~$n$ times a $1$ right of the
starting position of the head and for the remainders only~$0$s. Then
we execute the program, step by step. The machine might halt. This
is the case if after a finite number of steps, the program instructs
to stop the machine. If the tape consists of~$m$ times a~$1$ right
of the ending position of the head and for the rest~$0$s, we
call~$m$ the \emph{output}.

A partial function~$f\colon\N \to \N$ is~Turing-computable if there
is a program such that for every~$n$ where~$f$ is defined,
the Turing Machine loaded
with that program, halts on input~$n$ with output~$f(n)$ and
for every~$n$ where~$f$ is undefined, the Turing Machine loaded
with that program does not halt.

\subsection{Computable partial functions and sets}
Pick a preferred model of computation.
Programs (for Turing Machines),~$\lambda$-terms (in de~$\lambda$-Calculus)
and~$\mu$-recursive functions can all be enumerated effectively.
Let~$\varphi_e$ denote the~$e$th
\index{phi@$\varphi_e$}
computable partial function from~$\N^*$ to~$\N$,
where~$\N^*$ are the finite lists of natural numbers.
We define~$\varphi^k_e = \varphi_e \restriction \N^k \to \N$.
\index{phik@$\varphi^k_e$}

It is in general undecidable whether~$\varphi_e(x)\halts$. However,
if we limit ourselves to~$s$ steps of computation and define:
\[ \varphi_{e,s}(x) = \begin{cases}
                  \varphi_e(x) & \text{$\varphi_e(x)$ requires at
                                        most~$s$ steps of computation} \\
                  \divs & \text{$\varphi_e(x)$ did not stop
                                        after~$s$ steps of computation.}
                  \end{cases} \]
Then~$\varphi_{e,s}(x)\halts$ is decidable.

A set~$A\subseteq\N$ is called computable if its characteristic function
is computable. That is: there is an~$e$ such that~$\ind_A(x)=\varphi^1_e(x)$.

\subsubsection{Basic results from Recursion Theory}
In this thesis we will assume some familiarity with the basic
results of Recursion Theory, such as the following.\cite{odi1}
\begin{enumerate}
\item
Lists of natural numbers are effectively encodable with natural numbers. Thus
there exists a computable bijection~$\left< \  \right>\colon \N^* \to \N, 
(x_1, \ldots, x_m) \mapsto \left<x_1, \ldots, x_m\right>$.

\item
Arguments can be `hardcoded' in the program. That is:
for all~$n,m\in\N$, there is a
computable~$s^n_m\colon \N^{m+1} \to \N$,
such that for all~$x_1, \ldots, x_{m+n} \in \N$ we have:
\[
\varphi_{s^n_m(e,x_1, \ldots, x_m)}(x_{m+1},\ldots, x_{m+n}) 
        \simeq \varphi_e (x_1, \ldots, x_m, x_{m+1}, \ldots, x_{m+n} ).
\]
This result is
called~\keyword{Kleene's $s_{mn}$-Theorem}%
\index{smn-Theorem@$s_{mn}$-theorem}.\cite{klee38}
\item
Computation is computable. There are computable~$\operatorname{UM}$
and~$\operatorname{UM}'$ such that:\cite{tur37}\cite{klee38}
\begin{alignat*}{2}
\operatorname{UM}(p,x) & \simeq \varphi_p(x) &\qquad& (p,x\in \N)\\
\operatorname{UM}'(p,x,s). & \simeq \varphi_{p,s}(x) && (p,x,s\in\N)
\end{alignat*}
\end{enumerate}

\subsection{Recursively enumerable sets}
\begin{dfn}
Let~$\W_e = \{ n; \, n \in \N; \, \varphi_e(n) \halts \} = \dom \varphi^1_e$.
$\W_e$ is the
\keyword{$e$th recursively enumerable set}%
\index{enumerable!recursively!$e$th|see{$\W_e$}}%
\index{We@$\W_e$}.
Let~$\W_{e,s} = \{ n; \, n \in \N; \, \varphi^1_{e,s}(n) \halts \}$.
\end{dfn}
\begin{prop}[Kleene \cite{klee36}]
The following are equivalent
\begin{enumerate}
\item
$A = \W_e = \dom \varphi^1_e$ for some~$e$.

\item
$A = \im \varphi^1_e$ for some~$e$.

\item
$A$ has a computable enumeration. That is: there is an~$e$ such that:
\begin{itemize}
\item
if~$A$ is finite, then~$\varphi^1_e(n)\divs$ for~$n \geq \#A$
and~$A = \{ \varphi^1_e(0), \ldots, \varphi^1_e(n-1) \}$. 
\item
if~$A$ is infinite, then~$\varphi^1_e$ is total, injective
and~$A = \{ \varphi^1_e(0), \varphi^1_e(1), \ldots \}$.
\end{itemize}
\end{enumerate}
\end{prop}
\begin{proof}
\begin{description}
\item[1 $\Rightarrow$ 2]
Suppose~$A = \W_e$ for some~$e\in\N$.
There is a~$e'$ such that
\[ \varphi^1_{e'}(x) = \begin{cases}
                     x & \varphi^1_e(x) \halts \\
                     \divs & \varphi^1_e(x)\divs.
                     \end{cases} \]
Then~$x \in \im \varphi^1_{e'} \Leftrightarrow x \in \dom \varphi^1_e 
                \Leftrightarrow x \in A$.
\item[2 $\Rightarrow$ 3]
Suppose~$A = \im \varphi^1_e$ for some~$e$.
Consider the computable~$f$ given by:
\begin{algorithmic}[5]
\Function{$f$}{$n$}
\State \textbf{set} $A \gets \emptyset$
\State \textbf{set} $s \gets 0$
\For{$m$ \textbf{in} $\{ 0, \ldots, s\}$}
\If{$\varphi_e^s(m)\halts$ \textbf{and} $\varphi_e(m) \notin A$}
\State \textbf{append} $\varphi_e(m)$ \textbf{to} $A$
\If{$\#A = n+1$}
\State \Return $\varphi_e(m)$
\EndIf
\EndIf
\EndFor
\EndFunction
\end{algorithmic}
The function~$f$ returns~$x$ if and only if~$x = \varphi_{e,s}(m)$ for
some~$s,m\in\N$. Thus~$x \in A$. Furthermore, it is injective.
Finally, suppose~$x \in A$. Then~$\varphi_{e,t}(m)=x$ for some~$t,m\in\N$.
For some~$n$, in the execution of $f(n)$, the variable~$s$ will become bigger
than~$t$ and~$m$. $f(n)$ might return another number~$y$
if~$\varphi_{e,s}(m') \halts$ for a~$m' < m$. However, there are only
finitely many~$m'$ for which this can happen and for some~$n' > n$
we will have~$f(n')=x$. Thus~$f$ is the desired enumeration.
\item[3 $\Rightarrow$ 1]
Suppose~$\varphi_e$ is a computable enumeration of~$A$.
There is a~$e'\in\N$ such that~$\varphi_{e'}(n) = \mu t [\varphi_e(t) = n]$.
Then~$\varphi_{e'}(n)\halts$ if and only if~$\varphi_e$ enumerates~$n$ and
thus if and only if~$n \in A$. \qedhere
\end{description}
\end{proof}

\begin{prop}[Post \cite{post43}, Kleene \cite{klee36}
        and Mostowski \cite{mos47}]\label{posts-prop}
Given an~$A\subseteq \N$. If both~$A$ and~$\overline{A}$ are recursively
enumerable, then~$A$ is decidable.
\end{prop}
\begin{proof}
Let~$f\colon \N^3 \to \N$ be given by the following algorithm.
\begin{algorithmic}[5]
\Function{$f$}{$e,e',x$}
\State \textbf{set} $s \gets 0$
\Loop
\If{$x \in \W_{e,s}$}
\State\Return 1
\EndIf
\If{$x \in \W_{e',s}$}
\State\Return 0
\EndIf
\State \textbf{set} $s \gets s + 1$
\EndLoop
\EndFunction
\end{algorithmic}
There is a computable~$p\colon\N^2 \to \N$ such that
\[ \varphi_{p(e,e')}(x) = f(e,e',x). \]
Given~$e,e'\in\N$ such that~$\W_e \cap \W_{e'} = \emptyset$, then
\begin{align*}
        \varphi_{p(e,e')}(x) = 1 & \iff x \in \W_e \\
        \varphi_{p(e,e')}(x) = 0 & \iff x \in \W_{e'} \\
        \varphi_{p(e,e')}(x) \divs & \iff x \notin \W_e
                                \text{ and } x \notin \W_{e'}.
\end{align*}
Thus in partical, if~$\W_e = A$ and~$\W_{e'} = \overline{A}$,
then~$\ind_A = \varphi^1_{p(e,e')}$.
\end{proof}

\subsection{Oracles and Turing degrees}
\subsubsection{Oracles}
Given a~$A\subseteq\N$, let~$\varphi^{A}_e$
\index{phiA@$\varphi^A_e$}
denote the~$e$th ``computable'' partial function~$\N^* \to \N$ that is
allowed as extra computational step to ask questions about membership
of~$A$. One says that the ``computable function''  may consult an
oracle about~$A$.\index{oracle} In the case of Turing Machines, one
could add an extra tape with on that tape the characteristic function
of~$A$.\cite{tur39}
It is convenient and harmless to assume
that~$\varphi^{\emptyset}_e = \varphi_e$ for all~$e$.

We define~$\varphi^{k,A}_e = \varphi^A_e \restriction \N^k \to \N$
\index{phikA@$\varphi^{k,A}_e$}
and~$\W^A_e = \{ n; \, n \in \N; \,\varphi^{1,A}_e(n)\halts \}$.
\index{W@$\W^A_e$}

\subsubsection{Turing reductions}
If~$\ind_B = \varphi^{1,A}_e$ for some~$e$, one might say:
\begin{itemize}
\item
$B$ is recursive in~$A$.
\item
$B$ is computable with knowledge of~$A$.
\item
$B$ is
\keyword{Turing reducible}%
\index{reducible!Turing ($\leqT$)}
to~$A$ --- in symbols: $B \leqT A$.
\end{itemize}

$\leqT$ is transitive and reflexive. Thus we can study~$\PP(\N)/\leqT$,
the subsets of the natural numbers modulo Turing reducibility.
That is: the equivalence classes of the equivalence relation
\[ A \eqT B \quad \iff \quad A \leqT B \quad \text{and} \quad B \leqT A. \]
These equivalence classes are called the
\keyword{Turing degrees}%
\index{Turing degrees}.\cite{post48}

Given a computable set~$A$, then~$\ind_A = \varphi_e =
\varphi_e^\emptyset$.  Conversely, suppose~$\ind_A = \varphi_e^{1,B}$
for some~$A$ and a computable~$B$, then we can modify~$e$ to replace
every query to its oracle with an appriopriate computation and
thus~$\ind_A = \varphi^1_{e'}$ for some~$e'$. Thus the Turing degree
of~$\emptyset$ is the set of computable sets.

The Turing degree of~$\K$ is strictly above the degree of~$\emptyset$.
A degree is called enumerable, if there is a recursively enumerable~$A$ in
that degree.  There are no enumerable degrees above~$\K$:

\begin{prop}[Post \cite{post44}]
Given~$A \subseteq \N$.  If~$A$ is recursively enumerable, then~$A \leqT \K$. 
\end{prop}
\begin{proof}
Find~$e$ such that~$A=\W_e$. 
Let~$h\colon \N \to \N$ be a total computable function such that
$\varphi_{h(x)}(y) = \varphi_e(x)$ for all~$e,x,y\in\N$.  Then:
\[
h(x) \in \K \iff \varphi_{h(x)}(h(x)) \halts
\iff \varphi_e(x) \halts
\iff x \in \W_e.
\]
To complete the proof, let~$\varphi_a^\K (x) = \ind_{\{ x; \, h(x) \in \K \}}
= \ind_A$.
\end{proof}

This allows us to state Post's problem more succinctly:
\begin{prob}[Post \cite{post44}]
Is there a Turing degree between that of~$\emptyset$ and~$\K$?
\end{prob}

\subsection{Other reductions}
The Turing reductions we have seen so far,
are of a special kind:
\begin{dfn}[Post \cite{post44}]
$A$ is \keyword{$m$-reducible to~$B$}%
\index{reducible!$m$- ($\leqm$)} (in symbols: $A \leqm B$) if there is
a total computable~$f$ such that for all~$e$, we
have~$e \in A \iff f(e) \in B$.
\end{dfn}

\subsection{Fixed-Point theorem and second recursion theorem}

\subsubsection{Second recursion theorem}
The second recursion theorem states that we can write algorithms
that use their own codenumber.  More precisely:
\begin{thm}[Kleene \cite{klee38}]
There is a total recursive~$f$ such that for all~$e$ we have
\[ \varphi_{f(e)}(x) \simeq \varphi_e(f(e), x). \qquad (e,x \in \N) \]
\end{thm}
\begin{proof}
Let~$S$ be a total recursive function such
that~$\varphi_{S(x,y)}(z) \simeq \varphi_x(y,z)$;
$g$ such that~$\varphi_{g(e)}(z,x)=\varphi_e(S(z,z),x)$
and~$f(e) = S(g(e),g(e))$.
Then
\begin{align*}
\varphi_e(f(e),x) & \simeq \varphi_e(S(g(e),g(e)), x)  \\
                  & \simeq \varphi_{g(e)}(g(e),x) \\
                  & \simeq \varphi_{S(g(e),g(e))}(x)
                    \simeq \varphi_{f(e)}(x).\qedhere
\end{align*}
\end{proof}
It is easy and useful to strengthen the theorem.
\begin{cor}\label{sec-rec-cor}
For any~$n, m \in \N$,
there is a total function~$f_m\colon \N^{m+1} \to \N$
such that for all~$e,x_1, \ldots, x_n,y_1, \ldots, y_m \in \N$:
\[ \varphi_{f_m(e,y_1, \ldots, y_k)}(x_1, \ldots, x_n) \simeq
        \varphi_e(f_m(e,y_1, \ldots, y_m), x_1, \ldots, x_n). \]
\end{cor}
\begin{proof}
Given~$n,m\in\N$. There exist computable~$h_1$, $h_2$ and~$f_m$ such that
for all~$x_1,\ldots,x_n,y_1,\ldots,y_m,p,e\in\N$, we have
\begin{align*}
\varphi_{h_1(e)}(x_1,\ldots,x_n)
        & \simeq \varphi_e(\left<x_1,\ldots,x_n\right>)  \\
\varphi_{h_2(e,y_1,\ldots,y_m)}(p, \left<x_1,\ldots,x_n\right>)
        & \simeq \varphi_e(h_1(p), x_1, \ldots, x_n, y_1,\ldots,y_m) \\
f_m(e,y_1,\ldots,y_m) & = h_1(f(h_2(e,y_1,\ldots,y_m))).
\end{align*}
And thus:
\begin{align*}
&\varphi_{f_m(e,y_1,\ldots,y_m)}(x_1,\ldots,x_n) \\
 &\qquad\simeq\varphi_{h_1(f(h_2(e,y_1, \ldots, y_m)))}(x_1,\ldots,x_n) \\
 &\qquad\simeq\varphi_{f(h_2(e,y_1,\ldots,y_m))}(\left<x_1,\ldots,x_n\right>)\\
 &\qquad\simeq\varphi_{h_2(e,y_1,\ldots,y_m)}(f(h_2(e,y_1,\ldots,y_m)),
                                        \left<x_1,\ldots,x_n\right>)\\
 &\qquad\simeq\varphi_{e}(h_1(f(h_2(e,y_1,\ldots,y_m))),
                x_1,\ldots,x_n, y_1, \ldots, y_m),
\end{align*}
as desired.
\end{proof}

On first sight, the proof of the second recursion theorem seems magical.
We will study a closely related
theorem, called the Fixed-Point theorem, which is very similar to
a theorem with the same name in the $\lambda$-Calculus.

\subsubsection{Fixed-Point theorem}
\begin{thm}[Kleene \cite{klee36}, Turing \cite{tur37}, Curry \cite{cur42} and
                        Rosenbloom \cite{ros50}]
There is a~$\lambda$-term~$\mathcal{Y}$ such that for all $\lambda$-terms~$A$
\[ A (\mathcal{Y} A ) =_\beta \mathcal{Y} A. \]
\end{thm}
\begin{proof}
Define~$\omega$, $C$ and~$\mathcal{Y}$ as follows
\begin{align*}
\omega & = \lambda x.\, xx & C &= \lambda xyz.\, x(yz) & \mathcal{Y} & =
\lambda x.\, (Cx\omega)(Cx\omega).
\end{align*}
Then:
\begin{align*}
\mathcal{Y}A & =_\beta  (CA\omega)(CA\omega) \\
             & =_\beta A(\omega(CA\omega)) \\
             & =_\beta A((CA\omega)(CA\omega)) \\
             & =_\beta A(\mathcal{Y}A).
\end{align*}
\end{proof}

\begin{thm}[Kleene \cite{klee38}]
There is a total recursive~$f$ such that for all~$e$ such that~$\varphi_e$ is
total we have
\[ \varphi_{f(e)} \simeq \varphi_{\varphi_e(f(e))}.  \]      
\end{thm}
\begin{proof}
Define~$\omega$, C and~$f$ such that
\begin{align*}
\varphi_\omega(x) & \simeq \varphi_x(x)  &
\varphi_{\varphi_{C(x,y)}(z)} & \simeq \varphi_{\varphi_x(\varphi_y(z))} &
f(e) \simeq \varphi_{C(e, \omega)}(C(e, \omega))
\end{align*}
and~$C$ is computable and~$\varphi_{C(x,y)}$ is total. Then~$f$ is total and
\begin{align*}
\varphi_{f(e)} & \simeq \varphi_{\varphi_{C(e, \omega)}(C(e,\omega))} \\
& \simeq \varphi_{\varphi_{e}(\varphi_\omega(C(e, \omega)))} \\
& \simeq \varphi_{\varphi_{e}(\varphi_{C(e,\omega)}(C(e, \omega)))} \\
& \simeq \varphi_{\varphi_{e}(f(e))},
\end{align*}
as desired.
\end{proof}
There is a short proof of the Fixed-Point theorem
with the assumption of the second recursion theorem. And vice versa.
As an example, we will prove a strong version of the Fixed-Point
theorem that we will use often.
\begin{thm}\label{sfpt}
For every~$n\in\N$ there is a~$F\colon \N \to \N$ such that
for all~$e$, if~$\varphi_e$ is total, then for all
$z_1,\ldots,z_n\in\N$ we have:
\[
\varphi_{\varphi_{F(e)}(z_1, \ldots, z_n)} \simeq
\varphi_{\varphi_e(\varphi_{F(e)}(z_1,\ldots,z_n), z_1,\ldots,z_n)}.
\]
\end{thm}
\begin{proof}
There is an~$a\in\N$ and a computable~$F\colon \N \to \N$ such that
\begin{align*}
\varphi_a(p,x,z_1,\ldots,z_n) & \simeq
        \varphi_{\varphi_{e}(p,z_1,\ldots,z_n)}(x) \\
\varphi_{F_(e)}(z_1,\ldots,z_n) & =
        f_{n+1}(a,e,z_1,\ldots,z_n),
\end{align*}
where~$f_{n+1}$ is from the strong second recursion theorem
(Corollary~\ref{sec-rec-cor}). Thus:
\begin{align*}
\varphi_{\varphi_{F(e)}(z_1,\ldots,z_n)}(x)
 & \simeq\varphi_{f_{n+1}(a,e,z_1,\ldots,z_n)}(x) \\
 & \simeq\varphi_a(f_{n+1}(a,e,z_1,\ldots,z_n), x, e, z_1,\ldots,z_n) \\
 & \simeq\varphi_{\varphi_e(f_{n+1}(a,e,z_1,\ldots,z_n),z_1,\ldots,z_n )}(x) \\
 & \simeq\varphi_{\varphi_e(\varphi_{F(e)}(z_1,\ldots,z_n), z_1,\ldots,z_n)}(x).
                                                                        \qedhere
\end{align*}
\end{proof}
\subsection{Some conventions}\label{ss:conventions}
\begin{enumerate}
\item
Given a set~$A \subseteq \N$.
Where it does not confuse,
we write~$A(n)$ for~$\ind_A(n)$.
\item
Given a subset~$A\subseteq \N$, we write~$\overline{A}$ for~$\N - A$.
\item
Consider the sets~$\Pfin(\N), \N^2, \ldots$ and their standard bijections
with the natural numbers.  We call a map between any of these sets and/or~$\N$
computable if the lifted map between the natural numbers is computable.
\end{enumerate}

\clearpage
\section{Varieties of effective undecidability}

\subsection{Creative sets}
Post proved: if both~$A$ and~$\overline{A}$ are recursively enumerable,
then~$A$ is recursive.
Thus: if~$A$ is recursively enumerable but undecidable, every attempt to
enumerate the complement must go awry.

Consider~$\K$. Let~$\W_e \subseteq \overline{\K}$. Then~$e \notin \K$, because
if~$e \in \K$, then~$e \in \W_e \subseteq \overline{\K}$, which is absurd.
Thus also~$e \notin \W_e$. We know for all~$e$:
\[ \text{if \quad $\W_e \subseteq \overline{\K}$ \quad
        then \quad  $e \in \overline{\K} - \W_e$}. \]
Thus every attempt to enumerate the complement of~$\K$ forgets some
element we can effectively determine beforehand. This has inspired
the following definition.
\begin{dfn}[Post \cite{post44} and Myhill \cite{myh55}]
A recursively enumerable set~$A$ is
called~\keyword{creative}\index{creative} if there
exists a (total) computable~$f \colon \N \to \N$ such that for all~$e$
\[ \text{if \quad $\W_e \subseteq \overline{A}$ \quad
        then \quad  $f(e) \in \overline{A} - \W_e$}. \]
\end{dfn}

\begin{prop} \label{creative-not-simple}
The complement of a creative set contains an infinite recursively
enumerable subset.
\end{prop}
\begin{proof}
Suppose~$A$ is creative via~$f$.
Let~$w$ be the computable map $\Pfin (\N) \to \N$  such that~$\W_{w(A)} = A$ for
all~$A \in \Pfin (\N)$ and
\[ B_0 = \emptyset \qquad B_{n+1} = B_n \cup \{ f(w(B_n)) \} \qquad
        B = \bigcup_n B_n. \]
Then~$B$ is a recursively enumerable, infinite subset of~$\overline{A}$.
\end{proof}

Creative sets are certainly not a solution to Post's problem:

\begin{thm}[Myhill \cite{myh55}]\label{creative-complete}
A set is creative if and only if it is $m$-complete. 
\end{thm}
\begin{proof}
Suppose~$A$ is $m$-complete via~$f$; that is: $e \in \K \Leftrightarrow f(e) \in A$.
Let~$h$ be such that~$z \in \W_{h(e)} \Leftrightarrow f(z) \in \W_e$.
Suppose~$\W_e \subseteq \overline{A}$. Observe that:
\begin{enumerate}
\item 
If~$\W_e \subseteq \overline{A}$, then~$\W_{h(e)} \subseteq
\overline{\K}$, since~$z \in \W_{h(e)} \Leftrightarrow f(z) \in \W_e
\subseteq \overline{A} \Rightarrow z \in \overline{\K}$.

\item
If~$\W_e \subseteq \overline{\K}$ then~$e \in \overline{\K} - \W_e$,
for~$e\in \K$ implies~$e\in\W_e$ and
therefore~$\W_e \not \subseteq \overline{\K}$.
\end{enumerate}
And thus
\[ \W_e \subseteq \overline{A} \Rightarrow
   \W_{h(e)} \subseteq \overline{\K} \Rightarrow
   h(e) \in \overline{K} - \W_{h(e)} \Rightarrow
   f(h(e)) \in \overline{A} - \W_e, \]
which makes~$A$ creative via~$f \circ h$.

Conversely, suppose~$A$ is creative via~$f$.
The strong Fixed-Point theorem (Theorem \ref{sfpt})
ensures there is a computable~$g\colon\N \to \N$ such that
\[ \W_{g(z)} = \begin{cases}
               \{f(g(z)) \} & \text{if $z\in \K$} \\
               \emptyset & \text{otherwise.}
               \end{cases} \]
Then~$A$ is $m$-complete via~$f \circ g$:
\begin{itemize}
\item
Suppose~$z\in\K$. Then~$\W_{g(z)} = \{ f(g(z)) \}$. We must
have~$f(g(z)) \in A$,
since otherwise~$\W_{g(z)} \subseteq \overline{A}$ and
thus~$f(g(z)) \in \overline{A} - \W_{g(z)}$, which is absurd.
\item
Conversely, suppose~$z \in \overline{\K}$. Then~$\W_{g(z)} = \emptyset$.
Hence~$\W_{g(z)} \subseteq \overline{A}$.
And consequently~$f(g(z)) \in \overline{A} - \W_{g(z)} = \overline{A}$. \qedhere
\end{itemize}
\end{proof}

Later on we will look at weaker varieties of creativity. First we turn our
attention at a different class of undecidable sets.

\subsection{Simple sets}

Recall that decidable sets have recursively enumerable complements
and that the complement of a creative set is not recursively enumerable,
but contains an infinite recursively enumerable subset.

There are recursively enumerable sets that are undecidable because
their infinite complement does not contain an infinite recursively enumerable
set. These are called simple.

\begin{dfn}[Post \cite{post44}]
A recursively enumerable set~$A$ is
called~\keyword{simple}\index{simple} if~$\overline{A}$
is infinite but does not contain an infinite recursively enumerable subset.
\end{dfn}

Post invented and investigated simple sets and sets with even
stronger similar properties (hypersimple and hyperhypersimple sets)
as possible solutions to his problem.  Most simple sets that are
easy to construct tend to be~$T$-complete, as we will show at the
end of this subsection Post did look in the right direction, for
the solution of Friedberg, Muchnik and others after them are often
\emph{simple}.

\begin{exa}[Post \cite{post44}] \label{posts-simple}
To create a simple set~$S$ we must ensure that:
\begin{enumerate}
\item
$\overline{S}$ is infinite and
\item
for every~$e$, if~$\W_e$ is infinite, then~$\W_e \cap S \neq \emptyset$.
\end{enumerate}
We dovetail the computation of all recursively enumerable sets. For
every~$e$, we will put the first~$x > 2e$, for which we
compute~$x \in \W_e$, in~$S$. To wit:
\[ S = \{ x; \, \exists s, e [ x > 2e ; x\in\W_{e,s} \text{ and }
                 \nexists  y > 2e [ x \neq y \text{ and } y \in \W_{e,s} ]] \}. \]
\end{exa}

Simple sets are quite abundant: there is one in every recursive
enumerable degree.  It will take some work and new notions to show this.
However, this work will proof useful later on.

\subsubsection{Retraceable sets}

\begin{dfn}[Dekker and Myhill \cite{dek58}]
A principal enumeration~$a$ of a set~$A$
(that is: ~$A = \{ a_0 < a_1 < \ldots\}$ is
called
\keyword{retraceable}\index{retraceable} if there is a (total)
computable~$f\colon \N \to \N$ such that
\[ f(a_0) = a_0 \qquad f(a_{n+1}) = a_n. \]
Note that in general the enumeration~$a$ is not computable.
\end{dfn}

\begin{dfn}
A set~$A$ is called \keyword{immune}\index{immune}
if it is infinite, but does not contain
an infinite recursively enumerable subset.
\end{dfn}

\begin{rem}
A set~$A$ is simple if and only if it is recursively enumerable and its
complement is immune.
\end{rem}

\begin{prop}[Dekker and Myhill \cite{dek58}]
A retraceable set is either recursive or immune. 
\end{prop}
\begin{proof}
Suppose~$A$ is retraceable via~$f$ and not immune. Then there is
an infinite recursively enumerable~$B \subseteq A$.
We want to prove that~$A$ is recursive. That is, we want a decision method
whether~$x \in A$.
Find a~$b\in B$ such that~$x < b$. There is one, since~$B$ is infinite.
Then calculate~$f^1(b), f^2(b), \ldots, f^b(b)$. Note that these are all
the elements of~$A$ smaller than~$b$. Thus~$x$ is among them if and only
if~$x\in A$.
\end{proof}

\begin{cor} \label{recoret-simple}
An undecidable recursively enumerable set, which has a retraceable
complement, is simple.
\end{cor}

\subsubsection{Deficiency sets}

Let~$A$ be a undecidable recursively enumerable set enumerated
without repetition by a computable~$a\colon \N \to A$. Note that
if~$a$ is ascending ($a(0) < a(1) < \ldots$), then~$A$ is decidable.

\begin{dfn}[Dekker and Myhill \cite{dek58}]
The
\keyword{deficiency set}\index{deficiency set}
of a recursively enumerable set~$A$ with
respect to a computable enumeration without repetitions~$a\colon \N \to A$ is
\[ \mathcal{D}_a = \{ s; \, \exists t > s [a(t) < a(s)] \}. \]
\end{dfn}

\begin{rem}
The elements of~$\overline{\mathcal{D}_a}$ are called the
\keyword{true stages}\index{true stages} of
the enumeration.  They are interesting because:
\[ x \in A \iff x \in \{a(0), \ldots, a(s_0)\} \quad \text { where } \quad
        s_0 = \mu s \in \overline{\mathcal{D}_a} [a(s) > x].  \]
\end{rem}

This directly leads to the following proposition.

\begin{prop}[Dekker and Myhill \cite{dek58}]
For all enumerations~$a\colon \N \to A$ without repetitions of~$A$, we
have
\begin{enumerate}
\item
$A \eqT \mathcal{D}_a$

\item
$\overline{\mathcal{D}_a}$ is retraceable and thus (by corollary
\ref{recoret-simple}) $\mathcal{D}_a$ is simple.
\end{enumerate}
\end{prop}
\begin{proof}
$A \leqT \mathcal{D}_a$ by the preceding remark.  
Conversely, to decide whether $s \in \mathcal{D}_a$, given~$A$ we
check whether \[ A \cap \{0, \ldots, a(s) \} = \{ a(0), \ldots,
a(s) \}. \] This is precisely the case when~$s$ is a true stage and
thus~$s \notin \mathcal{D}_a$.

To show that~$\overline{\mathcal{D}_a}$ is retraceable, we have to compute the
element preceding a given~$s \in \overline{\mathcal{D}_a}$ (if it exists):
we pick the largest~$t < s$ such that
\[ \{ a(0), \ldots, a(s) \} \cap \{ 0, \ldots, a(t) \} = \{ a(0),
        \ldots, a(t) \}. \]
If there is no such~$t$, then there is no predecessor.
\end{proof}

\begin{cor}
There is a simple set in every undecidable recursively enumerable degree.
\end{cor}

\subsubsection{Effectively simple sets and Arslanov's criterium}

Although there are simple sets that fail to be~$T$-complete, any simple set for
which we know it is simple in an effective way \emph{is} $T$-complete.

\begin{dfn}[Smullyan \cite{smu64}]
A recursively enumerable set~$A$ is called
\keyword{effectively simple}\index{simple!effectively} if there
is a computable~$f\colon \N \to \N$ such that
\[ \text{if} \quad  \W_e \subseteq \overline{A} \quad \text{then} \quad
        \left| \W_e \right| \leq f(e). \]
\end{dfn}

\begin{exa}
Post's simple set (example \ref{posts-simple}) is effectively simple
via~$e\mapsto 2e+1$.
\end{exa}

We will first need to study the subtle completeness criterium of
Arslanov, which can be understood as a generalization of the Fixed Point
Theorem of Recursion Theory to any function of incomplete degree.

\begin{thm}[Martin \cite{mar66}, Lachlan \cite{lach68}
                        and Arslanov \cite{arsl81}] \label{arslanov}
A recursively enumerable set~$A$ is~$T$-complete if and only if
there exists a function~$f \leqT A$ without fixed points.
That is: for all~$e$, $\W_e \neq \W_{f(e)}$.
\end{thm}
\begin{proof}
Suppose~$A$ is~$T$-complete. Then we can compute~$\{ e; \, 0 \notin \W_e \}$
in~$A$.  Define~$f$ such that
\[ \W_{f(e)} = \begin{cases}
               \N & 0 \notin \W_e \\
               \emptyset & \text{otherwise.}
               \end{cases} \]
Then~$\W_{f(e)}$ and~$\W_e$ differ on~$0$ for every~$e$. Thus~$f$
is fixed-point free.

Conversely, suppose~$f \leqT A$ and~$f$ has no fixed points.
Let~$e$ be such that~$f = \varphi^{1,A}_e$. 
For every~$s$, the function~$\varphi^{1,A_s}_{e,s}$ is an approximation of~$f$,
but because it has fixed-points it differs from~$f$.

Let~$s$ be the modulus of~$\K$. That is:
\[ s_x = \begin{cases}
         \mu s [x \in \K_s] & x\in \K \\
         0 & \text{otherwise.} 
         \end{cases} \]
Note that if~$s \leqT A$ then~$\K \leqT A$.

Find a~$g$ using the Fixed Point Theorem such that
\[ \W_{g(x)} = \begin{cases}
             \W_{\varphi_{e, s_x}^{1,A_{s_x}}(g(x))} & x \in \K \\
             \emptyset & \text{otherwise.}
             \end{cases} \]
Since both~$f$ and~$g$ are total, we can compute~$\psi \leqT A$ such that
\[ f(g(x)) = \varphi_{e, \psi_x}^{1,A_{\psi_x}}(g(x)). \]
If~$x\in\K$ then~$g(x)$ is a fixed-point of~$\varphi^{1,A_{s_x}}_{e,s_x}$ and
thus~$s_x < \psi_x$. If~$x\notin \K$, then~$s_x = 0 < \psi_x$.
Consequently $x \in \K \Leftrightarrow x \in \K_{\psi_x}$.
Hence~$\K \leqT A$.
\end{proof}

\begin{prop}[Martin \cite{mar66}]
Every effectively simple set is~$T$-complete.
\end{prop}
\begin{proof}
Let~$A$ be effectively simple via~$f$. Define~$g \leqT A$ such that
\[ \W_{g(e)} = \{ \text{the first $f(e)+1$ elements of $\overline{A}$} \}. \]
Then~$g$ cannot have fixed-points: if~$\W_e = \W_{g(e)}$
then~$\W_e \subseteq \overline{A}$, however~$|\W_e| = f(e)+1$, which is absurd.
Thus~$A$ is~$T$-complete by the criterium of Arslanov.
\end{proof}

\subsection{Weakly and strongly effectively undecidable sets}

\subsubsection{Plain w.e.u.~and creative sets}

In the definitions of simple and creative sets we effectively thwart
properties of decidable sets. This is an indirect approach.
One could wonder what the status is of the following direct approach.

\begin{dfn}
A recursively enumerable set~$A$ is called
\keyword{(weakly) effectively undecidable}%
\index{w.e.u.}%
\index{effectively undecidable!weakly|see{w.e.u.}}
(\keyword{w.e.u.}) if there exists a
computable~$f\colon \N \to \N$ such that for all~$e$
\[ \text{if}\quad \varphi_e(f(e)) \halts \quad \text{then} \quad
        \varphi_e(f(e)) \neq A(f(e)). \]
\end{dfn}

The usual argument that~$\K$ is undecidable also shows
that~$\K$ is w.e.u. Actually, with a slight modification it
shows that every~$m$-complete set is w.e.u.

\begin{prop}
Every~$m$-complete set is w.e.u. 
\end{prop}
\begin{proof}
Suppose~$A$ is recursively enumerable and~$m$-complete via~$f$. That is:
$e \in \K \Leftrightarrow f(e) \in A$.
Let~$h$ be the computable function such that
\[ \varphi_{h(e)}(x) =  \begin{cases}
                        0 & \varphi_e(f(x)) = 0 \\
                        \divs & \text{otherwise.} 
                        \end{cases} \]
Suppose~$\varphi_e(f(h(e))) = A(f(h(e)))$. Then
\[ f(h(e)) \in A \Leftrightarrow \varphi_e(f(h(e))) \not\simeq 0
                 \Leftrightarrow \varphi_{h(e)}(h(e)) \divs
                 \Leftrightarrow h(e) \notin \K
                 \Leftrightarrow f(h(e)) \notin A, \]
which is absurd. Thus if~$\varphi_e(f(h(e))\halts$, then~$\varphi_e(f(h(e))
\neq A(f(h(e)))$. This shows~$A$ is w.e.u.~via~$f\circ h$.
\end{proof}

However, the w.e.u.~sets are not a new class.

\begin{prop}[Veldman]
Every w.e.u.~set is creative.
\end{prop}
\begin{proof}
Suppose~$A$ is w.e.u.~witnessed by~$f$. Let~$g\colon \N \to \N$ be
such that for all~$e$, if~$\W_e \subseteq \overline{A}$ then
\[ \varphi_{g(e)}(x) = \begin{cases}
                  0 & x \in \W_e \\
                  1 & x \in A \\
                  \divs & \text{otherwise.}
                  \end{cases} \]
\begin{itemize}
\item
Suppose~$\varphi_{g(e)}(f(g(e)))=0$. Then~$f(g(e)) \in \W_e$ by definition
of~$g$ and thus~$f(g(e)) \notin A$. However by w.e.u.~$\varphi_{g(e)}(f(g(e))) =
0 \neq 1 = A(f(g(e)))$. Contradiction.
\item
Suppose~$\varphi_{g(e)}(f(g(e)))=1$. Then~$f(g(e)) \in A$ by definition
of~$g$. However by w.e.u.~$\varphi_{g(e)}(f(g(e))) = 1 \neq 0 = A(f(g(e)))$.
Contradiction.
\end{itemize}
Apparently, as it is the only remaining possibility:
$\varphi_{g(e)}(f(g(e))) \divs$. And thus by definition of~$g$:
$f(g(e)) \notin \W_e$ and~$f(g(e)) \notin A$. Hence~$f(g(e)) \in
\overline{A} - \W_e$. This shows~$A$ is creative via~$f \circ g$.
\end{proof}

\subsubsection{$D$-w.e.u.~and $\W$-w.e.u.~sets}

The set constructed by Muchnik and Friedberg is not w.e.u.. We can
effectively determine at which point the set might differ from a
given decidable set, but we might be mistaken a recursively bounded
amount of times. More precisely:

\begin{dfn}
A recursively enumerable set~$A$ is
called~\keyword{$\W$-w.e.u.}%
\index{effectively undecidable!$\W$-weakly|see{$\W$-w.e.u.}}%
\index{W-w.e.u.@$\W$-w.e.u.}%
, if there
is a computable map~$f\colon \N \to \N ^2$ such that for all~$e$
\begin{enumerate}
\item
$| \W_{f(e)_1} | \leq f(e)_2$ and

\item
If $\forall z \in \W_{f(e)_1} [ \varphi_e (z) \halts ]$,
then $\exists z \in \W_{f(e)_1} [\varphi_e (z) \neq A(z)]$.
\end{enumerate}
\end{dfn}

Before we will investigate this class of sets, we will restrict our
attention to an easier subclass:

\begin{dfn}
A recursively enumerable set~$A$ is
called~\keyword{$D$-w.e.u.}%
\index{effectively undecidable!$D$-weakly|see{$D$-w.e.u.}}%
\index{D-w.e.u.@$D$-w.e.u.}%
, if there
is a computable map~$f\colon \N \to \Pfin(\N)$ such that for all~$e$
\[ \text{If}\quad\forall z \in f(e) [ \varphi_e (z) \halts ]
\quad\text{then}\quad\exists z \in f(e) [\varphi_e (z) \neq A(z)]. \]
\end{dfn}

One wonders whether there is a solution to Posts problem that is~$D$-w.e.u.
This is not the case:

\begin{thm}
Suppose~$a$ is a recursive enumeration without repetitions of
a~$D$-w.e.u.~set~$A$, then the deficiency set~$\mathcal{D}_a$ is
effectively simple.
Hence every~$D$-w.e.u.~set is~$T$-complete.
\end{thm}
\begin{proof}
Let~$A$ be~$D$-w.e.u.~via~$f$; $\alpha$ be such that~$A=\W_\alpha$
and~$g$ be the computable function such that
\[ z \in \W_{g(e)} \Leftrightarrow \exists y \in \W_e [a(y) > z\text{ and }
        z \notin \{ a(0), \ldots, a(y) \} ]. \]
Furthermore, let~$p$ be a function (see Proposition~\ref{posts-prop})
such that for all~$e_1,e_2 \in \N$,
if~$\W_{e_1} \cap \W_{e_2} = \emptyset$, then
\begin{align*}
\varphi_{p(e_1,e_2)}(z) = 1 & \ \Leftrightarrow \  z \in \W_{e_1} \\
\varphi_{p(e_1,e_2)}(z) = 0 & \ \Leftrightarrow \  z \in \W_{e_2} \\
\varphi_{p(e_1,e_2)}(z)\divs  & \ \Leftrightarrow \  
                                z \notin \W_{e_1} \cup \W_{e_2}.
\end{align*}
Suppose~$\W_e \subseteq \overline{\mathcal{D}_a}$.
Note that then~$\W_{g(e)} \subseteq \overline{A}$.
Let~$\xi = p(\alpha, g(e))$.

Assume~$\varphi_\xi (z) \halts$ for each~$z \in f(\xi)$.
Then there must be a~$z \in f(\xi)$ such that~$\varphi_\xi(z) \neq A(z)$.
Imagine~$\varphi_\xi(z)=1$. Then by
definition of~$p$ and~$\alpha$, we have~$z \in A$. However, by definition
of~$f$ also~$z\notin A$. Contradiction.
Apparently~$\varphi_\xi(z) = 0$. Thus~$f(z) \in A$, but then again by
definition of~$p$ we should have~$\varphi_\xi(z)=1$. Another contradiction.

Thus there is a~$z' \in f(\xi)$ such that~$\varphi_\xi (z')$ diverges.
Hence~$z' \notin A$ and~$z' \notin\W_{g(e)}$.
Thus for all~$y \in \W_e$,
we have~$a(y) < \max f(\xi)$ and consequently~$|\W_e| < \max f(\xi)$.
\end{proof}
\begin{cor}[Smullyan \cite{smu64}]
The deficiency set of~$\K$ is effectively simple.
\end{cor}

We proved that creative sets, $m$-complete sets and w.e.u.~sets
coincide.  Can we show that~$D$-w.e.u.~coincides with another notion
of completeness and creativity using similar proofs?  It turns out
we cannot use the same proofs for~$D$-w.e.u.~sets in general.
However, we can do it for a special subclass.

\subsubsection{$D$-s.e.u.,~$d$-complete and quasicreative sets}

\begin{dfn}
A recursively enumerable set~$A$ is
called~\keyword{$D$-strongly effectively undecidable}%
\index{effectively undecidable!$D$-strongly|see{$D$-s.e.u.}}%
\index{D-s.e.u.@$D$-s.e.u.}
(\keyword{$D$-s.e.u.}) via~$f\colon \N \to \Pfin (\N)$
if
\begin{enumerate}
\item
$A$ is $D$-w.e.u.~via~$f$ --- that is: for all~$e$ there is a~$z\in f(e)$
such that~$A (z) \neq \varphi_e(z)$ and furthermore
\item
$f(e)\subseteq \overline{A}$ for all~$e$. 
\end{enumerate}
\end{dfn}

\begin{dfn}[Shoenfield \cite{shoe57}]
A recursively enumerable set~$A\subseteq \N$ is
called~\keyword{quasicreative}%
\index{creative!quasi-}
if there is a computable function~$f\colon \N \to \Pfin(\N)$ such that for
all~$e$:
\[ \text{if}\quad \W_e \subseteq \overline{A} \quad \text{then} \quad
f(e) \subseteq \overline{A} \quad\text{and}\quad \exists z \in f(e) [
z \in \overline{A} - \W_e].  \]
\end{dfn}

\begin{dfn}[Jockusch \cite{jock66}]
A recursively enumerable set~$A$ is called~\keyword{disjunctive complete}%
\index{complete!disjunctive ($d$-)}
(\keyword{$d$-complete}) if there exists a computable~$f\colon \N \to \Pfin(\N)$
such that
\[ e \in \K \iff \exists z\in f(e) [z \in A]. \]
\end{dfn}

In \cite{shoe57} Shoenfield showed that~$d$-completeness and quasicreativeness
coincide.  Furthermore:

\begin{prop} \label{seu-equivs}
For a recursively enumerable~$A$, the following are equivalent:
\begin{enumerate}
\item
$A$ is~$D$-s.e.u.
\item
$A$ is quasicreative.
\item
$A$ is~$d$-complete.
\end{enumerate}
\end{prop}

\begin{proof}
\begin{description}
\item[$D$-s.e.u.~$\Rightarrow$ quasicreative]
Suppose~$A$ is~$D$-s.e.u.~via~$f$.
Let~$g\colon \N \to \N$ be
such that for all~$e$, if~$\W_e \subseteq \overline{A}$ then
\[ \varphi_{g(e)}(x) = \begin{cases}
                  0 & x \in \W_e \\
                  1 & x \in A \\
                  \divs & \text{otherwise.}
                  \end{cases} \]
We will prove that~$A$ is quasicreative via~$f \circ g$. Note
that for all~$e$ we have~$f(e)\subseteq\overline{A}$.
Assume that $\forall z\in f(e)[\varphi_{g(e)}(z)\halts]$.
Then there must be a~$z\in f(e)$ such that $\varphi_{g(e)}(z) \neq A(z)$,
but both
\begin{gather*}
\varphi_{g(e)}(z) = 1 \Rightarrow z \in A  \Rightarrow \varphi_{g(e)}(z) \neq 1 \\
\varphi_{g(e)}(z) = 0 \Rightarrow z \in \W_e 
                \Rightarrow z \notin A
        \Rightarrow \varphi_{g(e)}(z) \neq 0 
\end{gather*}
are absurd. Thus there is a~$z \in f(e)$ such that~$\varphi_{g(e)}(z)\divs$.
Then~$z\in \overline{A} - \W_e$ and we are done.

\item[quasicreative~$\Rightarrow$ $d$-complete, Shoenfield \cite{shoe57}]
Suppose~$A$ is quasicreative witnessed by~$f$. Let~$g$ be the computable
function such that
\[ \W_{g(e)} = \begin{cases}
               f(g(e)) & e \in \K \\
               \emptyset & \text{otherwise.}
               \end{cases} \]
Then~$A$ is~$d$-complete via~$f\circ g$:
\begin{itemize}
\item
If~$e \in \K$ then~$\W_{g(e)} = f(g(e))$. Also~$\exists z\in f(g(e))
[ z \in A]$, for otherwise~$\W_{g(e)}\subseteq\overline{A}$ and
then there is a~$z\in f(g(e))$ such that~$z\in \overline{A} - \W_{g(e)}$,
which is absurd since~$f(g(e)) = \W_{g(e)}$.

\item
Conversely, if $e \in \overline{\K}$ then $\W_{g(e)} = \emptyset$ and
$\W_{g(e)} \subseteq \overline{A}$ and thus there is a~$z\in f(g(e))$
such that~$z \in \overline{A} - \W_{g(e)} = \overline{A}$.
\end{itemize}

\item[$d$-complete $\Rightarrow$ $D$-s.e.u.]
Suppose~$A$ is~$d$-complete witnessed by~$f$. Let~$h$ be the computable function
such that
\[ \varphi_{h(e)}(x) = \begin{cases}
                       0 & \forall z\in f(x)[\varphi_e(z) = 0] \\
                       \divs & \text{otherwise.}
                       \end{cases} \]
Suppose~$\varphi_e(z) = A(z)$ for all~$z\in f(h(e))$. Then
\begin{align*}
        \exists z \in f(h(e))[ z\in A]
        & \Leftrightarrow \exists z \in f(h(e)) [\varphi_e(z) \neq 0] \\
        & \Leftrightarrow \varphi_{h(e)} (h(e)) \divs \\
        & \Leftrightarrow h(e) \notin \K \\
        & \Leftrightarrow \forall z \in f(h(e)) [z \notin A],
\end{align*}
which is absurd. Thus if for every~$z\in f(h(e))$ we know
$\varphi_e(z)\halts$, then there is a~$z\in f(h(e))$ such
that~$\varphi_e(z)\neq A(z)$ as desired. \qedhere
\end{description}
\end{proof}

\begin{prop}
The complement of a quasicreative set contains an infinite recursively
enumerable subset.
\end{prop}
\begin{proof}
Similar to the proof of proposition~\ref{creative-not-simple}.
\end{proof}

\begin{prop}[Shoenfield \cite{shoe57}]
There is a quasicreative set that is not creative.
\end{prop}

We will use a different and finer construction than Shoenfield's
to not only show that the notions~$d$-complete and~$m$-complete differ,
but also to show that several intermediate distinct completeness notions exist.

\subsubsection{$n$-$d$-complete sets}
\begin{dfn}
For a~$n\in\N$, the recursively enumerable set~$A$ is
called~$n$-$d$-complete via~$f\colon \N \to \Pfin(\N)$ if
\begin{enumerate}
\item
$A$ is $d$-complete witnessed by~$f$ --- that is: for all~$e$
\[ e \in \K \iff \exists z \in f(e) [z \in A]. \]
\item
$|f(e)| \leq n$ for all $e$.
\end{enumerate}
\end{dfn}
This completeness notion has a corresponding reduction:
\begin{dfn}
For a~$n\in\N$ and recursively enumerable sets~$A$ and~$B$, we
say~$A$ $n$-$d$-reduces to~$B$
(in symbols: \keyword{$A \leq_{n\textrm{-}d} B$})
\index{reducible!$n\textrm{-}d$- ($\leq_{n\textrm{-}d}$)}
via~$f\colon \N \to \Pfin(\N)$
if for all~$e\in\N$
\begin{enumerate}
\item
$|f(e)| \leq n$.
\item
$e \in A \iff \exists z \in f(e)[z \in B]$.
\end{enumerate}
\end{dfn}
\begin{rem}
\begin{enumerate}
\item
This reduction corresponds to the completeness:
$A$ is~$n$-$d$-complete if and only if every recursively
enumerable~$B$ $n$-$d$-reduces to~$A$.
\item
However, $\leq_{n\text{-}d}$ is not transitive. See
Corollary~\ref{n-d-not-transitive}.
\end{enumerate}
\end{rem}

\begin{prop}
For all~$i < j$ there is a set~$A$ that is~$j$-$d$-complete and
not~$i$-$d$-complete.\label{ijdcompl}
\end{prop}
\begin{proof}
We will define two recursively enumerable sets~$A$ and~$B$ such that
for all~$e\in\N$ the following hold:
\begin{description}
\item[$P_e$:]
$ e \in \K \iff je \in A \vee je+1\in A \vee \ldots \vee je + j -1\in A.$
\item[$N_e$:]
There is a~$p_e$ such that if~$\varphi_e(\left<e,p_e\right>) \halts$
and~$|D_{\varphi_e(\left<e,p_e\right>)}| \leq i$ then either
\begin{itemize}
\item
$\left< e, p_e\right> \in B$ but~$\forall z \in D_{
        \varphi_e(\left<e,p_e\right>)}[z \notin A]$ or
\item
$\left< e,p_e\right> \notin B$ but~$\exists z \in D_{
        \varphi_e(\left<e,p_e\right>)}[z \in A]$.
\end{itemize}
\end{description}
If~$P_e$ holds for all~$e\in\N$, $A$ is~$j$-$d$-complete.
If~$N_e$ holds, then~$B$ does not~$i$-$d$-reduce to~$A$ via~$\varphi_e$.
Thus if~$N_e$ holds for all~$e \in \N$, the set~$A$ is not~$i$-$d$-complete.

At some stage of the construction of~$A$ and~$B$, we notice
might that~$e \in \K$.
To meet~$P_e$, we will have to put one of~$\{je,\ldots,je+j-1\}$ in~$A$.
However, for every~$n\in\{0,\ldots,j-1\}$, there might be several~$e'$ such
that~$N_{e'}$ currently holds because~$je+n \notin A$.
For those~$N_{e'}$ that are `injured' because~$je+n$ is put in~$A$,
we will have to find a new~$p_{e'}$, which might involve new elements
we wish to keep out of~$A$.

To ensure that eventually all requirements are met, we will
allow a requirement to be injured only a finite number of times. To
that end, we use the so-called priority method.
The requirement $N_{e_1}$ is considered of higher priority than~$N_{e_2}$
if~$e_1 < e_2$. We will only allow a requirement to be injured in favor
of a requirement of higher priority.
Thus, when we need to choose~$je+n \in \{je, \ldots, je+j-1\}$ to put in~$A$,
we consider~$C_{je+n} = \{ e'; N_{e'}\text{ depends on } je+n \notin A\}$.
The priority by which we keep~$je+n$ out of~$A$ will be the
highest priority of~$C_{je+n}$.  That is~$\min C_{je+n}$.
We will put the~$je+n$ in~$A$ with the lowest priority.

Because~$i < j$, a requirement~$N_{e'}$ cannot depend on all~$\{je, \ldots,
je+j-1\}$ staying out of~$A$. Thus when~$e\in\K$ and~$N_{e'}$ has the highest
priority among the~$C_{je+n}$ it will not be injured.  
We say that the requirements that are injured, are injured in favor of~$e'$.
Note that a requirement will only be injured in favor of a requirement of
higher priority.

We will describe a program, that generates~$A$ and~$B$.
The program uses the following variables.
\begin{itemize}
\item
$C_e$ : A list of~$z$ such that  `$\varphi_z$ is not a $i$-$d$-reduction
                from~$B$ to~$A$' depends on~$e \notin A$
\item
$p_e$ : The number of times `$\varphi_e$ is not a $i$-$d$-reduction
                from~$B$ to~$A$' has injured.
\end{itemize}
Although there is an infinite amount of variables and the program
will perform operations on an infinite amount of them, we can
represent this with a finite structure: we do not store the values
of each of the variables, but the operations performed.
For instance, to represent the variables~$p_e$, we use a finite
list of pairs. Each pair~$\left<v ,e\right>$ means `we set all
variables~$p_i$ to~$v$ when~$\varphi_e(i)=1$'.  To get the current
value of~$p_i$, we simply replay the operations.

The program first initializes its variables and a list of conditions
to watch for: the watchlist.  Initially, the program watches
for~$e\in\K$ and~$\varphi(\left<e,0\right>)\halts$ for all~$e\in\N$.
After initialization, the program starts to compute on the conditions
in the watchlist in parallel.  When a condition is found to hold,
the specified action is executed and the condition is removed from
the watchlist.

Without further ado: the program to generate~$A$ and $B$.

\begin{algorithmic}[5]
\Initially
\State \textbf{add} $e\in\K$ \textbf{for all} $e \in \N
                \textbf{ to the watchlist}$
\State \textbf{add} $\varphi_e(\left< e, 0 \right>) \halts
                \textbf{ for all } e \in \N
                \textbf{ to the watchlist}$
\State \textbf{set} $p_e \gets 0$ \textbf{for all} $e \in \N$
\State \textbf{set} $C_e \gets \left<\right>$ \textbf{for all} $e \in \N$
\EndInitially

\When{$e\in\K$}
\If{there is a~$0 \leq n < j$ such that~$C_{je+n}=\left<\right>$}
\State \textbf{with that} $n$
\State \textbf{put} $je+n$ \textbf{in} $A$
\Else
\State \textbf{find} $n$ \textbf{such that}~$\min C_{je+n}$ is maximal in
        $\{ \min C_{je+m}; 0 \leq m < j \}$
\State \textbf{put} $je+n$ \textbf{in} $A$
\For{$z$ \textbf{in} $C_{je+n}$}
\State \textbf{set} $p_z \gets p_z + 1$
\State \textbf{add} $\varphi_z (\left<z, p_z\right>)\halts
                \textbf{ to the watchlist}$
\State \textbf{remove} $z$ \textbf{from} $C_a$ \textbf{for all} $a\in\N$
\EndFor
\EndIf
\EndWhen

\When{$\varphi_e(\left<e, p_e \right>) \halts$}
\State \textbf{set} $P \gets D_{\varphi_e(\left<e,p_e\right>)}$
\If{$|P| \leq i$ \textbf{and} $\forall z\in P[z \notin A]$}
\State \textbf{put} $\left< e, p_e \right>$ \textbf{in} $B$
\For{$z \in P$}
\State \textbf{append} $e$ \textbf{to} $C_z$
\EndFor
\EndIf
\EndWhen
\end{algorithmic}
We say that at a stage in the construction, the requirement~$N_e$ has settled
down if after that moment there will not be a requirement injured
in favor of~$e$.

All requirements~$N_e$ will settle. We will prove this by induction.
Suppose all~$e'<e$ will settle.
\begin{enumerate}
\item
Because all requirements of higher priority will settle, there will
be a stage after which~$e$ will not be injured anymore.        
Thus~$p_e$ will not change anymore.

\item
If for the final~$p_e$ it happens to be that
$\varphi_e(\left<e,p_e\right>)\divs$, the requirement~$N_e$ will
not pick elements which it wants to stay out of~$A$. And thus no other
requirement will be injured anymore in favor of~$N_e$.

In the other case, there will be a later stage
when it has been discovered that~$\varphi_e(\left<e,p_e\right>) \halts$
and acted upon. That is, the program might have put~$e\in C_z$ for several~$z$.
After that moment, $e$ will not be put in~$C_z$ anymore for any~$z$.

\item
In the latter case, there is a later stage such that
for all~$e'$ and~$0 \leq n < j$ such that~$e \in C_{je'+n}$,
either:
\begin{enumerate}
\item $e' \notin \K$ or
\item $e' \in \K$, which has been discovered and acted upon by the program.
This might involve injuring another requirement in favor of~$e$.
\end{enumerate}
After this stage, no requirement will be injured anymore in favor of~$e$.
Thus~$e$ settles down. \qedhere
\end{enumerate}
\end{proof}

\begin{cor}\label{n-d-not-transitive}
$\leq_{n\text{-}d}$ is not transitive.
\end{cor}
\begin{proof}
By the previous, we know there is a recursively enumerable~$A\subseteq\N$
such that~$A$ is~$2n$-$d$-complete, but not~$n$-$d$-complete.
In particular, there is a computable~$f\colon \N \to \N^{2n}$ such
that
\[ x \in \K \iff \exists i \in \{1,\ldots,2n\} [f(x)_i \in A]. \]
And now consider the~$B\subseteq \N$ such that
\begin{align*}
2x \in B & \iff \exists i \in \{1,\ldots,n\} [f(x)_i \in A] \\
2x+1 \in B & \iff \exists i \in \{n+1,\ldots,2n\} [f(x)_i \in A].
\end{align*}
Then~$\K \leq_{n\text-d} B \leq_{n\text-d} A$, but~$\K \nleq_{n\text-d} A$,
for otherwise~$A$ would be~$n$-$d$-complete.
\end{proof}

\begin{prop}
There is a set~$A$ that is~$d$-complete, but not~$n$-$d$-complete for all~$n$.
\end{prop}
\begin{proof}
We will use a variant of the program of the previous proof to
generate recursively enumerable~$A$ and~$B$ such that for all~$e,i\in \N$
the following holds
\begin{description}
\item[$P_e$:]
$e \in \K \iff \left<e,0\right> \in A \vee \ldots \vee \left<e,e \right> \in A.$
\item[$N^i_e$:]
There is a~$p^i_e$ such that if~$\varphi_e(\left<e,p^i_e\right>) \halts$
and~$|D_{\varphi_e(\left<e,p^i_e\right>)}| \leq i$ then either
\begin{itemize}
\item
$\left< e, p^i_e\right> \in B$ but~$\forall z \in D_{
        \varphi_e(\left<e,p^i_e\right>)}[z \notin A]$ or
\item
$\left< e,p^i_e\right> \notin B$ but~$\exists z \in D_{
        \varphi_e(\left<e,p^i_e\right>)}[z \in A]$.
\end{itemize}
\end{description}
The~$P_e$ combined assert that~$A$ is~$d$-complete. The
requirement~$N_e^i$ asserts~$\varphi_e$ is not a~$i$-$d$-reduction
from~$B$ to~$A$. Note that~$N^i_e$ implies~$N^j_e$ for~$j<i$, but
the converse does not hold.

Pick a computable bijection~$\pi\colon \N ^2 \to \N$.
The requirement~$N^i_e$
has higher priority than~$N^{i'}_{e'}$ (in symbols: $N^i_e \succ
N^{i'}_{e'}$) if~$\pi(i,e) < \pi(i',e')$.
\begin{algorithmic}[5]
\Initially
\State \textbf{add} $e\in\K$ \textbf{for all} $e \in \N
                \textbf{ to the watchlist}$
\State \textbf{add} $\varphi_e(\left< e, i, 0 \right>) \halts
                \textbf{ for all } e,i \in \N
                \textbf{ to the watchlist}$
\State \textbf{set} $p^i_e \gets 0$ \textbf{for all} $e \in \N$
\State \textbf{set} $C_e \gets \left<\right>$ \textbf{for all} $e \in \N$
\EndInitially

\When{$e\in\K$}
\If{there is a~$0 \leq n \leq e$ such that~$C_{\left<e,n\right>}=\left<\right>$}
\State \textbf{with that} $n$
\State \textbf{put} $\left<e,n\right>$ \textbf{in} $A$
\Else
\State \textbf{find} $n$ \textbf{such that}~$\min C_{\left<e,n\right>}$
        is maximal in
        $\{ \min C_{\left<e,m\right>}; 0 \leq m \leq e \}$
\State \textbf{put} $\left<e,n\right>$ \textbf{in} $A$
\For{$\left< z, i\right>$ \textbf{in} $C'_{\left<e,m\right>}$}
\State \textbf{set} $p^i_z \gets p^i_z + 1$
\State \textbf{add} $\varphi_z (\left<z, i, p^i_z\right>)\halts
                \textbf{ to the watchlist}$
\State \textbf{remove} $\left<z, i\right>$
                \textbf{from} $C'_a$ \textbf{for all} $a\in\N$
\State \textbf{remove} $\pi(z, i)$
                \textbf{from} $C_a$ \textbf{for all} $a\in\N$
\EndFor
\EndIf
\EndWhen

\When{$\varphi_e(\left<e, i, p^i_e \right>) \halts$}
\State \textbf{set} $P \gets D_{\varphi_e(\left<e,i,p^i_e\right>)}$
\If{$\forall z\in P[z \notin A]$ \textbf{and} $|P| \leq i$}
\State \textbf{put} $\left< e,i, p^i_e \right>$ \textbf{in} $B$
\For{$z \in P$}
                        \label{d-nd-odd-case2}
\State \textbf{append} $\left<e, i\right>$ \textbf{to} $C'_z$
\If{$\neg\exists e' \forall 0 \leq n \leq e' [\left< e', n \right> \in P]$}
\State \textbf{append} $\pi(e,i)$ \textbf{to} $C_z$
\EndIf
\EndFor
\EndIf
\EndWhen
\end{algorithmic}
The~$d$-completeness is obvious from the construction.
The~$i$-$d$-incompleteness is not.
Again, we say~$N^i_e$ settles at a stage if from that moment on
no other requirements will be injured in favor of~$N^i_e$.
Suppose all~$N^{i'}_{e'} \succ N^i_e$ eventualy settle.
\begin{enumerate}
\item
Because all requirements of higher priority will settle, there will
be a stage after which~$N^i_e$ will not be injured anymore in favor
of another requirement.

However, it may happen that~$N^i_e$
depends on all~$\{ \left<e', 0 \right>, \ldots , \left<e',e'\right> \}$
staying out of~$A$ for some~$e'$. If it is discovered that~$e' \in \K$,
the requirement~$N^i_e$ will be injured.

This may happen several times, but it will stop because of the following.
To be injured in this way, requires an~$e'$ such that at least
\begin{itemize}
\item $e' \in \K$, but this must not have been discovered and acted upon
at the moment~$\varphi_e(\left< e,i,p^i_e\right>)\halts$ is handled.
\item $e' \leq i$, because otherwise~$D_{\varphi_e(\left<e,i,p^i_e\right>)}$
cannot have less than or equal~$i$ elements and contain
all~$\{ \left<e', 0 \right>, \ldots , \left<e',e'\right> \}$.
\end{itemize}
Thus there are only a finite amount of oppurtunities for this kind of
injury. Hence, eventually, $N^i_e$ will not be injured anymore.
Thus~$p_e$ will not change anymore.
\end{enumerate}
The rest of the argument is the same as for the previous program.
\end{proof}

\subsubsection{$D$-w.e.u.~and hyper-simple sets}
Recall that~$D$-s.e.u.~sets are quasicreative; that quasicreative sets
contain an infinite recursively enumerable subset in their complement
and thus are not simple.
We will construct a simple~$D$-w.e.u.~set and thus show that~$D$-w.e.u.~and
$D$-s.e.u.~are different.

\begin{prop}
There is a simple~$D$-w.e.u.~set.
\end{prop}
\begin{proof}
We will modify the construction of Post's simple set (see
example~\ref{posts-simple}).
Let~$X_e$ be a computable partition of~$\N$ with~$|X_e|=e+2$. For instance, 
let~$X_e = \{ \frac{1}{2}(e+2)(e-1)-1 ,\ldots, \frac{1}{2}(e+3)(e+2)-2\}$.
For every~$e$
\begin{enumerate}
\item
Add the first~$x \in \W_e$
outside~$X_0 \cup \ldots \cup X_e$ to~$A$.
\item   
Add the first~$x \in X_e$ with~$\varphi_e(x)=0$ to~$A$.
\end{enumerate}
Clearly, $A$ is recursively enumerable.
Furthermore~$|X_e \cap \overline{A}| \geq 1$. By this and~(1),~$A$ is simple.
Let~$e$ be given. If there is an~$x \in X_e$ with~$\varphi_e(x) = 0$
then~$\ind_A$ will differ from~$\varphi_e$ on the first such~$x$ by~(2).
In the other case there is no~$x \in X_e$ with~$\varphi_e(x)=0$. Let~$x\in X_e$
be such that~$x \notin A$. Then~$\varphi_e(x)$ diverges or differs in value
with~$\ind_A$. Thus~$A$ is~$D$-w.e.u.~via~$e\mapsto X_e$.
\end{proof}

Although there is no recursively enumerable subset of the complement of
a~$D$-w.e.u.~set~$A$, we can recursively enumerate an infinite list of disjoint
finite sets, each of which intersects with the complement of our~$D$-w.e.u.~set.

\begin{dfn}[Post \cite{post44}]
\begin{enumerate}
\item
A computable~$h\colon \N \to \Pfin(\N)$ such that
if~$n \neq m$ then~$h(n) \cap h(m) = \empty$ is called a
\keyword{disjoint strong array}%
\index{disjoint strong array}.

\item
A set~$A$ \keyword{intersects with a disjoint strong array} $h$ if
for every~$n$ there is a~$z \in h(n)$ such that~$z\in A$.
\end{enumerate}
\end{dfn}

\begin{prop}
If~$A$ is~$D$-w.e.u., then there is a disjoint strong array
intersecting~$\overline{A}$.
\end{prop}
\begin{proof}
Suppose~$A$ is~$D$-w.e.u.~via~$f$.
Define~$\overline{w}$ to be a computable function such that for all finite
sets~$Y$:
\[
\varphi_{\overline{w}(Y)}(y) = \begin{cases}
                               0 & y \in Y \\
                               1 & y \notin Y.
                               \end{cases}
\]
Then~$f(\overline{w}(\emptyset))$ intersects~$\overline{A}$.
Thus our first try could be:
\[ X_0 = f(\overline{w}(\emptyset)) \qquad X_{n+1} =
f(\overline{w}(\bigcup_{i \leq n}X_i)). \]
This, however, does not work for it could be possible
that~$X_1 = X_0 - \overline{A}$. Secondly, $X_i$ is in
general not disjoint.

The following algorithm does generate a intersecting disjoint strong array:
\begin{algorithmic}[5]
\State \textbf{set} $G \gets \emptyset$ \Comment{current approximation of~$A$}
\State \textbf{set} $Y \gets \emptyset$ \Comment{already yielded}
\Loop
\State \textbf{set} $B \gets f(\overline{w}(G))$ \Comment{get the next batch}
\If {$B \cap Y \neq \emptyset$}
\State \textbf{set} $G \gets G - B$
\Else
\State \textbf{yield} $B$
\State \textbf{set} $Y \gets Y \cup B$
\State \textbf{set} $G \gets G \cup B$
\EndIf
\EndLoop
\end{algorithmic}
Note that in every iteration of the loop, by definition of~$f$,
there is a~$t\in B$ such that~$t \in A \cup G$ or~$t \notin A \cup G$.
If~$t \in A \cup G$, then$B \cap G \neq \emptyset$.
Note that~$G \subseteq Y$ and thus if~$Y \cap G = \emptyset$,
then with the previous~$B \cap \overline{A} \neq \emptyset$ and thus
we yield it.

Otherwise we remove the elements of~$B$ from~$G$.
Since~$G$ is finite, the algorithm cannot be stuck in this case.
Finally, since we check whether~$B \cap Y = \emptyset$ the yielded sets
are disjoint.
\end{proof}

\begin{dfn}
\begin{enumerate}
\item
A set~$A$ is called
\keyword{hyper-immune}
\index{immune!hyper-}
if it is infinite, but there is no
disjoint strong array intersection~$A$. 
\item
A set~$A$ is called
\keyword{hyper-simple}
\index{simple!hyper-}
if it is recursively enumerable
and has hyper-immune complement.
\end{enumerate}        
\end{dfn}

Thus $D$-s.e.u.~sets cannot be simple and~$D$-w.e.u.~sets cannot be hyper-simple.

In his 1944 article Post started his search for an intermediate set by
first considering simple sets, then hyper-simple sets and then even so
called hyperhyper-simple sets.

\subsubsection{$D$-w.e.u.~and $wtt$-complete sets}

We cannot generalize the reasoning of the proof of proposition~\ref{seu-equivs}
to~w.e.u.~sets. However, approaching the problem `from above' by specialing
$T$-completeness and the Arslanov's criterium is fruitful.

\begin{dfn}[Friedberg and Rogers \cite{frie59}]
A recursively enumerable set~$A$ is called
\keyword{$wtt$-complete}
\index{complete!$wtt$-}
if there is an~$e$ and a computable~$f\colon \N \to \N$ such
that~$\ind_\K = \varphi_e^{1,A}$ and the computation of~$\varphi_e^{1,A}(x)$
requires only queries to the oracle for elements below~$f(x)$.
\end{dfn}

\begin{prop}[Arslanov \cite{arsl81}]
A recursively enumerable set~$A \subseteq \N$ is~$wtt$-complete
if and only if there is a function~$f \leqwtt A$ without fixed-points.
\end{prop}
\begin{proof}
Follow the original proof (theorem~\ref{arslanov}) and note that
with the stronger assumptions we can replace~$T$-reduction
with~$wtt$-reductions.
\end{proof}

\begin{thm} \label{wtt-Dweu}
A set is~$D$-w.e.u.~if and only if it is~$wtt$-complete.
\end{thm}
\begin{proof}
Suppose~$A$ is~$D$-w.e.u.~witnessed by~$f$. Define~$g \leqwtt A$ such that
\[
\varphi_{g(e)}(x) = \begin{cases}
                    1 & x \in f(e) \text{ and } x \in A \\
                    0 & x \in f(e) \text{ and } x \notin A \\
                    \divs & \text{otherwise.}
                    \end{cases}
\]
Imagine~$g$ has a fixed-point~$e$, then for all~$z\in f(e)$ we
know~$\varphi_e(z) = \varphi_{g(e)}(z) = A(z)$. However, by~$D$-w.e.u., there
must be a~$z$ such that~$\varphi_e(z) \neq A(z)$. Contradiction.  
Thus~$g$ is fixed-point free and consequently by Arslanov's criteria
for~$wtt$-completeness, $A$ is $wtt$-complete.

Conversely, suppose~$A$ is~$wtt$-complete. Then there is an~$a$ and
a~$f\colon \N \to \N$ such that
$\K = \varphi_a^{1,A}$ and for all~$x$ the computation of~$\varphi_a^A(x)$ only
requires elements less than~$f(x)$ of~$A$.
Let~$e$ be given. There is a computable~$g\colon \N \to \N$
such that~$\varphi_a^{\varphi_e} \simeq \varphi_{g(e)}$. There is a
computable~$k\colon \N \to \N$ such that~$\K(k(e)) \not\simeq \varphi_e(k(e))$.
Thus
\[ \varphi_a^A(k(g(e))) \simeq  \K(k(g(e))) \not\simeq
        \varphi_{g(e)}(k(g(e))) \simeq
		\varphi_a^{\varphi_e}(k(g(e))). \]
The computation of~$\varphi_a^A(k(g(e)))$ only requires elements less
than~$f(k(g(e)))$ from~$A$, thus~$\varphi_e$ and~$A$ must differ on the
first~$f(k(g(e)))$ places. That is~$A$ is~$D$-w.e.u.~witnessed
by~$x\mapsto \{0, 1, \ldots, f(k(g(x)))\}$.
\end{proof}

\subsubsection{$wtt$-complete and weakly quasicreative sets}
We saw~$d$-completeness corresponds to quasicreativeness.  There
is a weakening of quasicreativeness to which~$wtt$-completeness
corresponds.
\begin{dfn}[Kanovich \cite{kano70}]
A recursively enumerable set~$A$ is called
\keyword{weakly quasicreative}
\index{creative!weakly quasi-}
if there exists a computable~$f \colon \N \to \Pfin(\N)$ such that for
all~$e$, if~$\W_e \subseteq \overline{A}$ then there is a~$z\in f(e)$
such that~$z \in \overline{A} - \W_e$.
\end{dfn}
Note that the difference with definition of quasicreative is
the absence of the requirement that~$f(e) \subseteq \overline{A}$
for all~$e\in\N$.

\begin{thm}[Kanovich \cite{kano70}]
A set is weakly quasicreative if and only if it is~$wtt$-complete.
\end{thm}
\begin{proof}
Suppose~$A$ is weakly quasicreative witnessed by~$f$. Define~$g
\leqwtt A$ such that~$\W_{g(e)} = f(e) \cap \overline{A}$.
Imagine $g$ has a fixed-point~$e$. Then~$\W_e = \W_{g(e)} =
\overline{A} \cap f(e)$, but there must be a~$z \in f(e)$ such
that~$z \in \overline{A} - \W_e = \emptyset$. Contradiction.
Thus~$g$ has no fixed-points and consequently~$A$ is~$wtt$-complete
by Arslanov's criterium.

Conversely, suppose~$A$ is~$wtt$-complete. Then by
Theorem~\ref{wtt-Dweu} it is~$D$-w.e.u.~via some~$g$.
Let~$h\colon \N \to \N$ be a computable function such that
\[ \varphi_{h(e)} = \begin{cases}
                    0 & x \in \W_e \\
                    1 & x \in A \\
                    \divs & \text{otherwise}.
                    \end{cases} \]
Suppose~$\W_e \subseteq \overline{A}$.
There is a~$z \in g(h(e))$ such that $\varphi_{h(e)}(z) \not\simeq A(z)$.
If~$z \in A$ then~$\varphi_{h(e)}(z) = 1 = A(z)$ and also
if~$z \in \W_e$, then $\varphi_{h(e)}(z) = 0 = A(z)$.
Thus~$z \in \W_e - \overline{A}$, which shows~$A$ is weakly quasicreative
via~$g\circ h$.
\end{proof}

\clearpage
\section{Solutions to Post's problem}
We have considered various notions of effective undecidability and have
proven that all of them imply~$T$-completeness.  That is, except
for~$\W$-w.e.u, but every undecidable set is~$\W$-w.e.u.

A solution to Post's problem is a reursively enumerable set~$A$
such that~$\emptyset \leT A \leT \K$.  We will review various
constructions of such sets and try to discover how far we can
stretch our effective knowledge of their undecidability.

\subsection{Two recursively enumerable sets of incomparable degree}
\label{incomp-degrees}
\begin{thm}[Friedberg \cite{frie57}]
There exist recursively enumerable sets~$A$ and~$B$ such
that~$A \not\leqT B$ and~$B \not\leqT A$.
\end{thm}

Since for any~$X$ and all computable~$C$, we know $C \leqT X$,
both~$A$ and~$B$ are undecidable. Furthermore, since for all complete~$K$
and recursively enumerable~$X$ we know~$X \leq T$, both~$A$ and~$B$ are
not $T$-complete. Thus both~$A$ and~$B$ are solutions to Post's problem.

\subsubsection{Sketch of the construction}
For every~$e$ we will try to find~$a_e$ and~$b_e$ for which we can ensure
\begin{enumerate}
\item $\varphi^B_e(a_e) \divs$ or $\varphi^B_e(a_e) \neq A(a_e)$ and
\item $\varphi^A_e(b_e) \divs$ or $\varphi^A_e(b_e) \neq B(b_e)$.
\end{enumerate}
For instance, when we discover that~$\varphi_{37}^{B_s}(a_{37} = 0$,
we will put~$a_{37}$ into~$A$. Here $B_s$ is~$B$ as far as it has
been defined at the current stage in the construction. The computation
of~$\varphi_{37}^{B_s}$ depends on a finite part of~$B$.  We must make sure
that~$B_s$ does not differ with~$B$ on that part, for otherwise
the~$\varphi_{37}^{B_s}$ might not equal~$\varphi_{37}^B$.

There might be a~$b_e$ in that finite part of~$B$ for which
we have yet to handle~$\varphi_{e}^{A_s}(b_e) \halts$. It might be the
case that we would like to put~$b_e \in B$, but by doing that would
alter the result of~$\varphi_{37}^{B}$.
We could change~$b_e$ to
a new value outside of the finite part of~$B$ on
which~$\varphi_{37}^{B_s}$ depends.
However, if we always do this, we might move~$b_e$ an infinite amount
of times.

To solve this, we will use priorities: ($x \prec y$ means $y$ has higher
priority than~$x$)
\[ a_0 \succ b_0 \succ a_1 \succ b_1 \succ \ldots \]
If~$b_e$ has lower priority than~$b_{37}$, then we will move~$b_e$.
Otherwise, we will move~$b_{37}$.

\subsubsection{Algorithm}
In the construciton we will use the same event-based setup we used in
Proposition~\ref{ijdcompl}. We will use the following variables:
\begin{enumerate}
\item
$a_e$: The element for which we try to ensure
$\varphi^B_e(a_e) \not\simeq A(a_e)$.
\item
$C^a_e$: $1$ when we handled~$\varphi^B_e(a_e)$ and~$0$ otherwise.
\item
$U^a_e$: the initial number of elements of~$A$ on
which~$\varphi^A_e(b_e)$ depends if~$C^b_e=1$ and~$-1$ otherwise.
\end{enumerate}
$b_e$, $C^b_e$ and~$U^a_e$ are used similarly with the r\^oles
of~$a$/$A$ and~$b$/$B$ swapped.
If we write~$A$ or~$B$ in the algorithm, we refer to~$A$ and~$B$ as far
as they have been defined at that moment.

\begin{algorithmic}[5]
\Initially
\State \textbf{add `} $\varphi_e^A(b_e)\halts$ and $C^b_e=0$
                \textbf{' for all} $e \in \N
                \textbf{ to the watchlist}$
\State \textbf{add `} $\varphi_e^B(a_e)\halts$ and $C^a_e=0$
                \textbf{' for all} $e \in \N
                \textbf{ to the watchlist}$
\State \textbf{set} $a_e \gets e$ \textbf{for all} $e \in \N$
\State \textbf{set} $b_e \gets e$ \textbf{for all} $e \in \N$
\State \textbf{set} $U^a_e \gets -1$ \textbf{for all} $e \in \N$
\State \textbf{set} $U^b_e \gets -1$ \textbf{for all} $e \in \N$
\State \textbf{set} $C^a_e \gets 0$ \textbf{for all} $e \in \N$
\State \textbf{set} $C^b_e \gets 0$ \textbf{for all} $e \in \N$
\EndInitially

\When{$\varphi_e^A(b_e)\halts$ and $C^b_e=0$}
\State \textbf{set} $C^b_e \gets 1$.
\State \textbf{find} $u$ \textbf{such that}
        $\varphi_e^A(b_e)$ depends on at most the first~$u$ elements of~$A$.
\State \textbf{set} $U^a_e \gets u$.
\For{$e' > e$ \textbf{such that} $a_{e'} \leq u$ and $C^a_{e'}=0$}
\State \textbf{set} $a\gets \mu a [A(a)=0 \text{ and } \max_x U^a_x \leq a
                                  \text{ and } \forall e'' [a_{e''}\neq a]]$.
\State \textbf{remove `} $\varphi_{e'}^B(a_{e'})\halts$ and $C^a_{e'}=0$
                \textbf{' from the watchlist}
\State \textbf{set} $a_{e'}\gets a$
\State \textbf{add `} $\varphi_{e'}^B(a_{e'})\halts$ and $C^a_{e'}=0$
                \textbf{' to the watchlist}
\EndFor
\For{$e' > e$ \textbf{such that} $U^b_{e'} \geq b_e $ and $C^a_{e'}=1$}
\State \textbf{set} $C^a_{e'}\gets 0$
\State \textbf{set} $U^b_{e'}\gets -1$
\State \textbf{set} $b\gets \mu a [A(a)=0 \text{ and } \max_x U^a_x \leq a
                                  \text{ and } \forall e'' [a_{e''}\neq a]]$.
\State \textbf{set} $a_{e'}\gets a$
\State \textbf{add `} $\varphi_{e'}^B(a_{e'})\halts$ and $C^a_{e'}=0$
                \textbf{' to the watchlist}
\EndFor
\EndWhen

\When{$\varphi_e^B(a_e)\halts$ and $C^a_e=0$}
\State Similarly, with~$A/a$ swapped with~$B/b$.
\EndWhen
\end{algorithmic}

\subsubsection{Analysis of effective undecidability}
$a_0$ is never moved. $b_0$ is only moved in favor of~$a_0$.
Thus at most once. $a_1$ is only moved in favor of~$b_0$. $b_0$ will
cause at most two moves of~$a_1$, thus~$a_1$ is moved at most twice.
Et cetera.

In general, let~$M_{a_e}$ denote the number of times~$a_e$ is moved.
Then
\[ M_{a_e} = \sum_{b_{e'} \succ a_e} M_{b_{e'}}+1 \qquad 
   M_{b_e} = \sum_{a_{e'} \succ b_e} M_{a_{e'}}+1.\]

\begin{prop}
$M_{a_n} = F_{2(n+1)}-1$ and~$M_{b_n} = F_{2(n+1)+1}-1$,
where~$F_n$ is the~$n$th Fibonacci number. That is:
\[ F_n = \begin{cases}
         0 & n = 0 \\
         1 & n = 1 \\
         F_{n-2} + F_{n-1} & n > 1.
         \end{cases} \]
\end{prop}
\begin{proof}
By induction.
\begin{description}
\item[Base]
\begin{align*}
M_{a_0} & = 0 = F_2 - 1 = F_{2(0+1)} - 1 \\
M_{b_0} & = M_{a_0}+1 = 1 = F_3 - 1 = F_{2(0+1)+1} - 1.
\end{align*}
\item[Step]
Suppose~$M_{a_n} = F_{2(n+1)}-1$ and~$M_{b_n} = F_{2(n+1)+1}-1$.
Then
\begin{align*}
M_{a_{n+1}} & = M_{a_n} + M_{b_n} + 1 \\
            & = F_{2(n+1)}-1 + F_{2(n+1)+1} -1 + 1 \\
            & = F_{2((n+1)+1)}-1 \\
M_{b_{n+1}} & = M_{b_n} + M_{a_{n+1}} + 1 \\
            & = F_{2(n+1)+1}-1 + F_{2(n+2)} -1 + 1 \\
            & = F_{2(n+2)+1} - 1. \qedhere
\end{align*}
\end{description}
\end{proof}
For the final~$a_e$ we know~$\varphi^B_e(a_e) \neq A(e)$. If~$e$ does
not use the oracle, then $\varphi_e(a_e) \neq A(a_e)$.

\subsection{A $n^2$-$\W$-w.e.u solution}

\begin{dfn}
A recursively enumerable set is called~$f$-$\W$-w.e.u if
there is a computable~$g\colon \N \to \N$ such that for all~$e$
\begin{enumerate}
\item $|\W_{g(e)}| \leq f(e)$
\item If~$\forall z \in \W_{g(e)}[\varphi_e(z)\halts]$,
        then~$\exists z \in \W_{g(e)}[\varphi_e(z) \neq A(z)]$.
\end{enumerate}
\end{dfn}
Thus the set~$A$ constructed in Subsection~\ref{incomp-degrees}
is~$(F_{2(n+1)}-1)$-$\W$-w.e.u.
We can modify the construction a bit to get a sharper result.
In the original construction the priorityorder of the
requirements is:
\[ a_0 \succ b_0 \succ a_1 \succ b_1 \succ \ldots \]
Given~$A_0,A_1,\ldots\in\N$. If we rearrange the priorities as follows.
\begin{align*}
a_0 \succ \ldots \succ a_{A_0-1}  \succ b_0 & 
\succ a_{A_0} \succ \ldots \succ a_{A_0+A_1-1} \succ b_1 \\
        & \succ a_{A_0+A_1} \succ \ldots
\end{align*}
Then~$M_{a_e} = M_{a_{e+1}} = \cdots = M_{a_{e+A_e-1}}$.
Define
\[ \alpha_n = M_{a_{A_0 + \cdots + A_{n-1}}} \quad \text{and}
                \quad \beta_n = M_{b_n}. \]
It is easy to verify that for~$n>0$:
\begin{align*}
        \alpha_0 &= 0 & \beta_0 &= A_0 \\
        \alpha_n &= \alpha_{n-1} + \beta_{n-1}+1 &
        \beta_n &= \beta_{n-1} + A_n (\alpha_n + 1).
\end{align*}
And thus
\begin{align*}
\alpha_n & = \alpha_{n-1} + 1 + \sum_{i = 0}^{n-1} A_i (\alpha_i + 1)  \\
         & \leq \alpha_{n-1} + 1 + \sum_{i = 0}^{n-1} A_1(\alpha_{n-1}+1) \\
         & =  (A_0 + \ldots + A_{n-1} + 1)(\alpha_{n-1}+1).
\end{align*}
And consequently:
\begin{prop}
Given any order-preserving
computable~$f\colon \N \to \N$ such that for every~$N \in \N$,
there exists a~$M \in \N$ such that
\[ f(M+N) \geq (N+M+1)(f(N) + 1). \]
There is an~$A\subseteq \N$
that is~$f$-$\W$-w.e.u and~$\emptyset \leT A \leT \K$.
\end{prop}
\begin{proof}
It is sufficient
to find~$A_0,A_1,\ldots\in\N$ such that for any~$n$ we
have
\[ f(A_0+\cdots+A_{n-1}) \geq \alpha_n. \]
By assumption on~$f$ here is a~$N$ such that
\begin{align*}
  f(N+0) & \geq (N+1)(f(0)+1)
\intertext{and thus}
             & \geq (N+1) \\
             & \geq (N+1)(\alpha_0+1) \\
             & \geq \alpha_1.
\end{align*}
Pick this~$N$ as~$A_0$.

Suppose we already found suitable~$A_0,\ldots,A_{n-1}\in\N$.
Then by assumption on~$f$ there is a~$N$ such that:
\begin{align*}
        & f(A_1+\cdots + A_{n-1}+N) \\
        &\qquad \geq (A_1+\cdots + A_{n-1}+N+1)(f(A_1+\cdots+A_{n-1})+1)
\intertext{and thus}
        &\qquad \geq (A_1+\cdots+A_{n-1}+N+1)(\alpha_n+1) \\
        &\qquad \geq \alpha_{n+1}.
\end{align*}
Thus use~$N$ for~$A_n$.
\end{proof}
\begin{cor}
There is a~$n^2$-$\W$-w.e.u solution to Post's problem.
That is: there is a recursively enumerable~$A$
such that~$A$ is~$n^2$-$\W$-w.e.u
and~$\emptyset \leT A \leT \K$.
\end{cor}

\clearpage
\section{Conclusion}
First we investigated some new direct notions of effective undecidability.
They turned out to be equivalent to previously investigated notions
of completeness and creativity.\footnote{
The notions in paranthesis were not covered in this thesis.}
\begin{center}
\begin{tabular}{ccc}
completeness & creativity & effective undecidability \\
\hline
$wtt$-complete & weakly quasicreative & $D$-w.e.u. \\
$d$-complete & quasicreative & $D$-s.e.u. \\
& \vdots & \\
$n$-$d$-complete & ($n$-creative) & ($n$-$D$-s.e.u.) \\
& \vdots & \\
$2$-$d$-complete & ($2$-creative) & ($2$-$D$-s.e.u.) \\
$m$-complete & creative & w.e.u.
\end{tabular}
\end{center}
This leads to the following question.
\begin{oprob}
Is there a notion of effective undecidability that is equivalent
to~$T$-completeness?
\end{oprob}

Then we investigated the effective undecidability of the
existing constructions of sets that are solutions to Post's problem.
We saw we could find a~$n^2$-$\W$-w.e.u solution. Can we do better?
\begin{oprob}
Is there a recursively enumerable~$A$ such that $\emptyset \leT A \leT \K$ and
$A$ is~$f$-$\W$-w.e.u.~for a bounded computable~$f\colon \N \to \N$?
\end{oprob}

\clearpage
\appendix
\section{Bibliography and index}
\bibliography{main}{}
\bibliographystyle{amsalpha}

\clearpage

\printindex

\end{document}